\documentclass[11pt]{article}

\usepackage{amsmath,amssymb,amsthm,upref}
\usepackage[ansinew]{inputenc}
\usepackage{dsfont}
\usepackage{mathrsfs}
\usepackage[all]{xy}
\usepackage{html}
\usepackage{fullpage}
\usepackage{paralist}

\newcommand{\cdcop}[1]{\mrm{C}^+({#1})}
\newcommand{\cdc}[1]{\mrm{C}_\ast({#1})}
\newcommand{\cdcp}[1]{\mrm{C}_+({#1})}

\newcommand{\Dl}{\ensuremath{\Delta} }
\newcommand{\simp}{\Delta^\comp}

\newcommand{\ds}{\displaystyle}

\newcommand{\mc}[1]{\mathcal{#1}}
\newcommand{\mrm}[1]{\mathrm{#1}}
\newcommand{\mbf}[1]{\mathbf{#1}}

\newcommand{\comp}{{\scriptscriptstyle \ensuremath{\circ}}}

\theoremstyle{theorem}
\newtheorem{thm}{\textsc{Theorem}}[section]

\newtheorem{lema}[thm]{\textsc{Lemma}}
\newtheorem{prop}[thm]{\textsc{Proposition}}
\newtheorem{cor}[thm]{\textsc{Corollary}}

\theoremstyle{definition}
\newtheorem{defi}[thm]{\textsc{Definition}}
\newtheorem{ej}[thm]{\textsc{Example}}

\newtheorem{obs}[thm]{\textbf{Remark}}

\newenvironment{num}{\medskip
\refstepcounter{thm}\noindent {\bf (\thethm)}}{\vspace{0.2ex}\par}

\setdefaultleftmargin{0em}{2.5em}{.5em}{.5em}{.5em}{.5em}

\begin{document}

\title{Simplicial descent categories}

\author{
\begin{tabular}{c} Beatriz Rodr\'iguez Gonz\'alez\footnote{The author was supported by the research projects `ERC Starting Grant TGASS', `Geometr\'{i}a Algebraica, Sistemas Diferenciales y Singularidades' FQM-218, MTM2007-66929, and by FEDER}\mbox{ {} }\footnote{email:rgbea@icmat.es}\\
\begin{small}ICMAT, Consejo Superior de Investigaciones Cient\'{i}ficas\end{small}\end{tabular}}

\date{\empty}

\maketitle

\begin{abstract}
Let $\mc{D}$ be a category and $\mrm{E}$ a class of morphisms in $\mc{D}$. In this paper we study the question of how to transfer homotopic structure from the category of simplicial
objects in $\mc{D}$, $\simp\mc{D}$, to $\mc{D}$  through a `good' functor $\mbf{s}:\simp\mc{D}\rightarrow\mc{D}$, which we call simple functor. For instance, the Bousfield-Kan
homotopy colimit in a Quillen simplicial model category is a good simple functor. As a remarkable example outside the setting of Quillen models we include Deligne simple of mixed Hodge complexes. We prove
here that the simple functor induces an equivalence on the corresponding localized categories. We also describe a natural structure of Brown category of cofibrant objects on $\simp\mc{D}$.
We use these facts to produce cofiber sequences on the localized category of $\mc{D}$ by $\mrm{E}$, which give rise to a natural Verdier triangulated structure in the stable case.
\begin{flushright}\begin{tabular}{r}
\begin{small}{\textit{2000 Mathematics Subject Classification:} 14F35 (Primary); 18G30, 18D99 (Secondary)}\end{small}\\
\begin{small}{keywords: Homotopy colimit; Brown category of cofibrant objects; $\Dl$-closed class} \end{small}
\end{tabular}
\end{flushright}
\end{abstract}

\section*{Introduction}
\mbox{}\indent Since the beginnings of homotopy theory the usage of simplicial techniques has been extremely fruitful. In the present paper we study how
to induce homotopic structure on a category $\mc{D}$
endowed with a class $\mrm{E}$ of morphisms using a `simple' functor from the category of simplicial objects in $\mc{D}$, $\simp\mc{D}$, to $\mc{D}$. Two classical precedents are the
following: the \textit{geometric realization} $\simp Top\rightarrow Top$, and
the `\textit{total complex}' of a double chain complex, which may be seen as a simple functor $\simp \cdc{\mc{A}} \rightarrow \cdc{\mc{A}}$.

In the cosimplicial setting, P. Deligne introduces in \cite{DeIII} a simple functor for cosimplicial mixed Hodge complexes, and uses it as a tool to define a mixed Hodge
structure over the cohomology of singular varieties. Given a singular variety $S$, P. Deligne constructs a mixed Hodge structure on $H^\ast(S)$ from $\mbf{s}(K^\bullet)$, which is the
simple mixed Hodge complex of a cosimplicial mixed Hodge complex $K^\bullet$ associated with a smooth hyperresolution $X_\bullet$ of $S$.  It turns out that the resulting mixed Hodge
structure on $H^\ast(S)$ is independent of the hyperresolution $X_\bullet$ chosen for $S$. However, a priori, it could depend on the choice of Deligne simple $\mbf{s}(\cdot)$.\\[0.1cm]
\indent A natural question is to find out the properties that make Deligne's construction $\mbf{s}$ a `good'
simple functor, and if this construction is unique or not. In this paper we introduce the notion of simplicial descent category as an answer to this question.\\
Consider a category $\mc{D}$ endowed with a saturated class $\mrm{E}$ of morphisms, called \textit{weak equivalences}, such that both $\mc{D}$ and $\mrm{E}$ are closed by finite
coproducts. A \textit{simplicial descent structure} on $(\mc{D},\mrm{E})$ is a triple
$(\mbf{s},\mu,\lambda)$ satisfying the following axioms:
\begin{quote}
(S1) \textsc{Coproducts}: $\mbf{s}:\simp\mc{D}\rightarrow\mc{D}$ is a functor, called the `simple functor', which commutes with finite coproducts up to weak equivalence.\\
(S2) \textsc{Eilenberg-Zilber}: $\mu$ is a zigzag of natural weak equivalences between the iterated simple of a bisimplicial object and the simple of its diagonal.\\
(S3) \textsc{Normalization}: $\lambda$ is a zigzag of natural weak equivalences between an object $A$ of $\mc{D}$ and the simple of the constant simplicial object induced by $A$.\\
(S4) \textsc{Exactness}: The simple of a degreewise weak equivalence is a weak equivalence.\\
(S5) \textsc{Homotopy}: If $e$ is a simplicial homotopy equivalence in $\simp\mc{D}$, then $\mbf{s}(e)$ is in $\mrm{E}$.
\end{quote}
A \textit{simplicial descent category} is a pair $(\mc{D},\mrm{E})$ as before, endowed with a simplicial descent structure. \textit{Cosimplicial descent categories} are defined dually, involving a simple $\mbf{s}:\Dl\mc{D}\rightarrow\mc{D}$.
This notion is strongly inspired by the \textit{cubical homological descent categories} introduced by F. Guill\'en and V. Navarro in \cite{GN}.\\[0.15cm]
\indent It holds that Deligne's construction provides a simple functor verifying our axioms.\\[0.1cm]
\noindent {\bf Theorem \ref{HodgeDescenso}}
\textit{Let $\mc{H}dg$ be the category of mixed Hodge complexes, and consider the class $\mrm{E}_{\mc{H}dg}$ of quasi-isomorphisms in $\mc{H}dg$. Then Deligne simple functor $\mbf{s}_{\mc{H}dg}:\Dl \mc{H}dg\rightarrow\mc{H}dg$ endows $(\mc{H}dg,{\mrm{E}_{\mc{H}dg}})$ with a structure of cosimplicial descent category.}\\[0.15cm]
\indent On the other hand, homotopy colimits in Quillen simplicial model categories are also examples of simple functors verifying our axioms. More concretely, we prove the\\[0.1cm]
\noindent {\bf Theorem \ref{Modelos}}
\textit{Let $\mc{M}$ be a simplicial model category with weak equivalences $\mc{W}$, and denote by $\mc{M}_c$ and $\mc{M}_f$ its corresponding subcategories of cofibrant and fibrant objects.
Then, $(\mc{M}_c,\mc{W},\mrm{hocolim})$ is a simplicial descent category. Dually, $(\mc{M}_f,\mc{W},\mrm{holim})$ is a cosimplicial descent category.}\\[0.15cm]
In addition we see in corollary \ref{ModelosCompl} that, under mild conditions on $\mc{M}$, this simplicial descent structure may be extended to all $\mc{M}$, by taking a `corrected' Bousfield-Kan homotopy colimit as simple functor.\\[0.1cm]
\indent The following facts highlight some remarkable differences between simplicial descent categories and Quillen model categories.
\vspace{0.1cm}
\begin{compactitem}[-]
\item The localized category of a simplicial descent category may not have small hom's (remark \ref{NoQ}).
\item Simplicial descent structures are inherited by diagram categories (proposition \ref{DescensoFuntores}).
\item A simplicial descent structure on $(\mc{D},\mrm{E})$ is unique, up to unique isomorphism, on the localizations (corollary \ref{unicSimple}).
\item Simplicial descent structures are closed by homotopical equivalence (proposition \ref{EquivHomot}).
\end{compactitem}

\vspace{0.15cm}

In the present paper we show that the previous axioms ensure good homotopic properties on $(\mc{D},\mrm{E})$. This is done in two steps. First, we study the homotopic structure of $(\simp\mc{D},\mc{S}=\mbf{s}^{-1}\mrm{E})$ and second, we prove that the simple functor transfers this structure. Indeed, it induces an equivalence between the corresponding localized categories of $\simp\mc{D}$ and $\mc{D}$.\\
\indent Independently, V. Voevodsky introduces in \cite{Vo} the notion of $(\Dl,\amalg_{< \infty})$-closed class of $\simp\mc{C}$, for a general category $\mc{C}$. This notion and the one of simplicial descent category are closely related.

\noindent {\bf Proposition \ref{DlClosed}}\\
\vspace{-0.4cm}
\textit{\begin{compactitem}
\item[\textbf{\textit{i}.}] Consider a category $\mc{C}$ with finite coproducts, and a saturated class $\mc{W}$ of morphisms in $\simp\mc{C}$. Then, the following are equivalent.
\begin{compactitem}
\item[1.] $\mc{W}$ is $(\Dl,\amalg_{< \infty})$-closed.
\item[2.] $(\simp\mc{C},\mc{W},\mrm{D}:\simp\simp\mc{C}\rightarrow\simp\mc{C})$ is a simplicial descent category.
\end{compactitem}
\item[\textbf{\textit{ii}.}] Given a simplicial descent category $(\mc{D},\mrm{E},\mbf{s})$, then
$\mc{S}=\mbf{s}^{-1}\mrm{E}$ is $(\Dl,\amalg_{< \infty})$-closed.
\end{compactitem}}

\vspace{0.15cm}

A $(\Dl,\amalg_{< \infty})$-closed class $\mc{W}$ of $\simp\mc{C}$ provides a natural structure of cofibrant objects in the sense of K. Brown \cite{Br}. In this case, the simplicial homotopic structure of $\simp\mc{C}$ is compatible with $\mc{W}$, in such a way that the natural cofiber sequences in $(\simp\mc{C})[\mc{W}^{-1}]$ satisfy the usual properties. Although not stated explicitly,
the proof of the following result is contained in \cite{Vo}.\\[0.1cm]
\noindent {\bf Proposition \ref{VovWald}}\\
\textit{
\textbf{\textit{i.}} Let $\mc{C}$ be a category with finite coproducts and $\mc{W}$ a $(\Dl,\amalg_{< \infty})$-closed class of $\simp\mc{C}$. Then
$(\mc{C},\mc{W},Cof)$ is a category of cofibrant objects, where $Cof=\{$termwise coprojections$\}$.\\
\textbf{\textit{ii.}} If in addition $\mc{C}$ is pointed, then a pair $\left( X\rightarrow Y \rightarrow Z  \, , \, Z\rightarrow Z\sqcup \Lambda(X) \right)$ is a simplicial distinguished triangle if and only if it is a cofibration sequence in the sense of Brown.}\\[0.15cm]
\indent These results reveal that the previous axioms guarantee good homotopic properties on $\simp\mc{D}$. But our aim is to work on $\mc{D}$, not on $\simp\mc{D}$. The key result that
makes possible to transfer structure from $\simp\mc{D}$ to $\mc{D}$ is the following.\\[0.1cm]
\noindent {\bf Theorem \ref{equivCat}}\mbox{}\\
\textit{\textbf{\textit{i.}} The simple functor $\mbf{s}:(\simp\mc{D})[(\simp\mrm{E})^{-1}]\rightarrow \mc{D}[\mrm{E}^{-1}]$ is left adjoint to
$c:\mc{D}[\mrm{E}^{-1}]\rightarrow (\simp\mc{D})[(\simp\mrm{E})^{-1}]$.\\
\textbf{\textit{ii.}} The pair $\mbf{s}:\simp\mc{D}\rightleftarrows \mc{D}:c$ is a homotopical equivalence between
$(\simp\mc{D},\mc{S})$ and $(\mc{D},\mrm{E})$. In particular, $\mbf{s}:\simp\mc{D}[\mc{S}^{-1}]\rightarrow \mc{D}[\mrm{E}^{-1}]$ is an equivalence
of categories.}\\[0.15cm]
The previous results imply the following corollary.\\[0.1cm]
\noindent {\bf Corollary \ref{EstructuraTriangulada}}\mbox{}\\
\textit{In the pointed case, distinguished triangles in $\mc{D}[\mrm{E}^{-1}]$ satisfy the usual `non-stable' axioms for triangulated categories, and they are
natural with respect to diagram categories.\\
Consequently, in the stable case, $\mc{D}[\mrm{E}^{-1}]$ has a structure of Verdier triangulated category which is natural with respect to diagram categories.}\\[0.2cm]
\indent There are some further interesting questions concerning the theory presented here. For instance, theorem \ref{equivCat} implies that
$\mbf{s}:\simp\mc{D}\rightarrow \mc{D}$ is a homotopy colimit. Indeed, in \cite{R2}
we obtain homotopy colimits $I\mc{D}\rightarrow\mc{D}$ for diagrams of finite shape $I$ by combining $\mbf{s}$ with the simplicial replacement.
In addition, in case $(\mc{D},\mrm{E})$ is closed by small coproducts, to have a simplicial descent structure on $(\mc{D},\mrm{E})$ turns out to be equivalent
to have realizable homotopy colimits of arbitrary shape on $(\mc{D},\mrm{E})$.\\[0.2cm]
\indent The organization of the paper goes as follows. We introduce the basic definitions and properties of simplicial descent categories in the first section. In the second one we
describe some concrete examples, as chain complexes, simplicial sets or mixed Hodge complexes. In section 3 we prove that the Bousfield-Kan homotopy colimit in a
Quillen simplicial model category is a simple functor in our sense. In section 4 we describe the relation between simplicial descent categories, $\Dl$-closed classes and Brown structures of cofibrant objects. Finally, in the last section we prove that the simple functor is an
homotopical equivalence. We use this result to deduce the existence and good properties of the induced cofiber sequences in $\mc{D}$.\\[0.25cm]
\indent I wish to express my deep gratitude to my thesis advisors L. Narv\'aez Macarro and V. Navarro Aznar, who suggested me this topic and gave me their helpful advice and dedication.

I am very grateful to an anonymous referee who observed that, under our axioms, $\simp\mc{D}$ supports a Brown structure of cofibrant objects, and who proposed to give an alternative
proof of corollary \ref{EstructuraTriangulada} using this fact, instead the direct proof given in a previous version of this paper. I would like also to thank the other referees
for useful comments and for stating the question whether the simple functor induces an equivalence on localized categories or not, which is now solved in theorem \ref{equivCat}.

\section*{Notations and preliminaries}

\begin{num}
Denote by $\Dl$ the \textit{simplicial category}, with objects the ordered sets $[n]=\{0,\ldots,n\}$, $n\geq 0$, and morphisms the order preserving maps. The \textit{face maps} $d^i:[n-1]\rightarrow [n]$ are characterized by $d^i([n-1])=[n]-\{i\}$, and the degeneracy maps $s^j:[n+1]\rightarrow [n]$ are the surjective monotone maps with $s^j(j)=s^j(j+1)$.
They satisfy the well-known \textit{simplicial identities}, and generate all maps in $\Dl$ (see, for instance, \cite{May}).

By $\simp\mc{D}$ (resp. $\simp\simp\mc{D}$) we mean the category of \textit{simplicial} (resp. \textit{bisimplicial}) \textit{objects} in a fixed category $\mc{D}$.
The \textit{diagonal functor} $\mrm{D}:\simp\simp\mc{D}\rightarrow\simp\mc{D}$ is given by $\mrm{D}(\{Z_{n,m}\}_{n,m\geq 0})=\{Z_{n,n}\}_{n\geq 0}$.\\
The \textit{constant simplicial object} defined by $A\in\mc{D}$ will be denoted by $c(A)$. In this way we obtain the constant functor $c:\mc{D}\rightarrow \simp\mc{D}$, which is fully faithful. When understood, we will denote $c(A)$ by $A$.

Dually, $\Dl\mc{D}=(\simp\mc{D}^\comp)^\comp$ is the category of \textit{cosimplicial objects} in $\mc{D}$.
\end{num}
\begin{num} If $\mc{D}$ has (finite) coproducts, there is a natural action of simplicial (finite) sets on
$\simp\mc{D}$, given by
$$(K\boxtimes X)_n=\coprod_{K_n}X_n$$
Recall that $\Dl[k]$ is the simplicial finite set with $\Dl[k]_n=Hom_{\Dl}([n],[k])$. In particular $\Dl[0]=\ast$, and $\Dl[1]$ plays the role of `unit interval' in $\simp Sets$. Given $X\in\simp\mc{D}$, the maps $d^0,d^1:[0]\rightarrow [1]$ induce $d_0^X,d_1^X:X\rightarrow X\boxtimes \Dl[1]$. Simplicial homotopies and simplicial homotopy equivalences are defined in $\simp\mc{D}$ as usual, using the simplicial cylinder $X\boxtimes\Dl[1]$ of $X$.
\end{num}
\begin{num}
We call a class of morphisms $\mrm{E}$ of a category $\mc{D}$ \textit{saturated} if $\mrm{E}=\gamma^{-1}(\mrm{isomorphisms})$, where $\gamma:\mc{D}\rightarrow\mc{D}[\mrm{E}^{-1}]$ is the localization functor. In particular, $\mrm{E}$ contains all isomorphisms of $\mc{D}$, is closed by retracts and satisfies the 2-out-of-3 property.\\
Given another category $\mc{C}$, $Fun(\mc{C},\mc{D})$ denotes the category of functors from $\mc{C}$ to $\mc{D}$. If $I$ is a small category, we also write $Fun(I,\mc{D})=I\mc{D}$. We denote by $\mc{C}\mrm{E}$, or by $\mrm{E}$ if $\mc{C}$ is understood, the class of morphisms $\tau:F\rightarrow G$ in $Fun(\mc{C}, \mc{D})$ such that $\tau_c\in\mrm{E}$ for all $c\in\mc{C}$.
\end{num}

\begin{num} By a \textit{relative category} we mean a pair $(\mc{D},\mrm{E})$ formed by a category $\mc{D}$ and a class of morphisms
$\mrm{E}$ of $\mc{D}$, which is assumed to be saturated. The morphisms of $\mrm{E}$ will be called \textit{weak equivalences}. We say that a relative
category $(\mc{D},\mrm{E})$ is \textit{closed by finite coproducts} if $\mc{D}$ has an initial object $0$ and both $\mc{D}$ and $\mrm{E}$ are closed by
finite coproducts.\\
Note that if $(\mc{D},\mrm{E})$ is a relative pair closed by finite coproducts then $\mc{D}[\mrm{E}^{-1}]$ is closed by finite coproducts, and they are preserved by
$\gamma:\mc{D}\rightarrow \mc{D}[\mrm{E}^{-1}]$.
\end{num}

%%%%%%%%%%%%%%%%%%%%%%%%%%%%%%%%%%%%%%%%%%%%%%%%%%%%%%%%%%%%%%%%%%%%%%%%%%%%%%%%%%%%%%%%%%%%%%%%%%%%%%%%%%%%%%%%%%%%%%%%%%%%%%%%%%
%%%%%%%%%%%%%%%%%%%%%%%%%%%%%%%%%%%%%%%%%%%%%%%%%%%%%%%%%%%%%%%%%%%%%%%%%%%%%%%%%%%%%%%%%%%%%%%%%%%%%%%%%%%%%%%%%%%%%%%%%%%%%%%%%%
%%%%%%%%%%%%%%%%%                   CAPITULO
%%%%%%%%%%%%%%%%%%%%%%%%%%%%%%%%%%%%%%%%%%%%%%%%%%%%%%%%%%%%%%%%%%%%%%%%%%%%%%%%%%%%%%%%%%%%%%%%%%%%%%%%%%%%%%%%%%%%%%%%%%%%%%%%%%

\section{(Co)simplicial Descent Categories}

\begin{defi}\label{defiCDS} Let $(\mc{D},\mrm{E})$ be a relative pair closed by finite coproducts.
A \textit{simplicial descent structure} on $(\mc{D},\mrm{E})$ is a triple $(\mbf{s},\mu,\lambda)$
satisfying the following five axioms.
\begin{compactitem}
\item[$\mathbf{(S1)}$] $\mbf{s}:\simp\mc{D}\rightarrow \mc{D}$ is a functor, called the \textit{simple functor}, which commutes with finite coproducts
up to weak equivalence. That is, the canonical morphism
$\mbf{s}(X) \amalg \mbf{s}(Y)\rightarrow \mbf{s}(X\amalg Y)$ is in $\mrm{E}$ for all $X$, $Y$ in $\simp\mc{D}$.
\item[$\mathbf{(S2)}$] $\mu:\mbf{s}\comp\mrm{D}\dashrightarrow\mbf{s}\comp\mathbf{s}$ is a zigzag of natural weak equivalences. If
$Z\in\simp\simp\mc{D}$, recall that $\mbf{s}\mrm{D}(Z)$ is the simple
of the diagonal of $Z$. On the other hand
$\mbf{s}\mbf{s}(Z):=\mbf{s}(n\rightarrow\mbf{s}(m\rightarrow
Z_{n,m}))$ is the iterated simple of $Z$.
\item[$\mathbf{(S3)}$] $\lambda:\mbf{s}\comp c\dashrightarrow
1_{\mc{D}}$ is a zigzag of natural weak equivalences, which is assumed to be compatible with $\mu$ in the sense of (\ref{compatibLambdaMu}) below.
\item[$\mathbf{(S4)}$] If $f:X\rightarrow Y$ is a morphism in
$\simp\mc{D}$ such that $f_n \in \mrm{E}$ for all $n$,
then $\mbf{s}(f)\in\mrm{E}$.
\item[$\mathbf{(S5)}$] The image under the simple functor of the map $d_0^A:A\rightarrow A\boxtimes \Dl[1]$ is a weak equivalence for each object $A$ of $\mc{D}$.
\end{compactitem}
\end{defi}
\begin{num}\label{compatibLambdaMu} Compatibility between $\lambda$ and $\mu$.\\
Given $X\in\simp\mc{D}$, denote by $X\times\Dl$, $\Dl\times X$ the bisimplicial objects with
$(X\times\Dl)_{n,m}=X_n$ and $(\Dl\times X)_{n,m}=X_m$. Note that $\mbf{s}\mbf{s}(X\times\Dl)=
\mbf{s}(n\rightarrow \mbf{s}c(X_n))$ and
$\mbf{s}\mbf{s}(\Dl\times X)= \mbf{s}c\mbf{s}(X)$. The
compositions
\begin{equation}\label{compatibLambdaMuEquac}\xymatrix@M=4pt@H=4pt@R=7pt@C=25pt{
 \mbf{s}(X) \ar[r]^-{\mu_{\Dl\times X}} &  \mbf{s}c\mbf{s}(X)\ar[r]^-{\lambda_{\mbf{s}(X)}} & \mbf{s}(X) &
 \mbf{s}(X) \ar[r]^-{\mu_{X\times \Dl}} &  \mbf{s}\mbf{s}c(X)\ar[r]^-{\mbf{s}(\lambda_{X})} &
 \mbf{s}(X)}\end{equation}
give rise to isomorphisms of $\mbf{s}$ in
$Fun(\simp\mc{D},\mc{D})[\mrm{E}^{-1}]$. Then, $\lambda$ is said to be \textit{compatible} with $\mu$ if the above isomorphisms are the identity in $Fun(\simp\mc{D},\mc{D})[\mrm{E}^{-1}]$.
\end{num}

\vspace{0.2cm}

\begin{defi}
A \textit{simplicial descent category} is a relative pair $(\mc{D},\mrm{E})$ closed by finite coproducts and endowed with a simplicial descent structure $(\mbf{s},\mu,\lambda)$.\\
Dually, a \textit{cosimplicial descent structure} on $(\mc{D},\mrm{E})$ is a triple $(\mbf{s}:\Dl\mc{D}\rightarrow\mc{D},\mu,\lambda)$ such that
$(\mbf{s}^\comp,\mu^\comp,\lambda^\comp)$ is a simplicial descent structure on $(\mc{D}^\comp,\mrm{E}^\comp)$. \textit{Cosimplicial descent categories} are defined analogously.
\end{defi}

To shorten the notations, we will also write  $(\mc{D},\mrm{E},\mbf{s})$ for a simplicial descent category $(\mc{D},\mrm{E})$ endowed with a simplicial descent structure
$(\mbf{s},\mu,\lambda)$.

\begin{obs}\mbox{}\\
\textsc{I.} If $(\mc{D},\mrm{E})$ admits a simplicial descent structure $(\mbf{s},\mu,\lambda)$, then $(\mbf{s},\mu,\lambda)$ is uniquely determined on $\mc{D}[\mrm{E}^{-1}]$ up to unique isomorphism. This is proved later in corollary \ref{unicSimple}.\\
\textsc{II.} Since $\mrm{E}$ is saturated, (S3) implies that the canonical map $0\rightarrow \mbf{s}(0)$ is a weak equivalence.\\
\textsc{III.} All results concerning simplicial descent categories are dualized to cosimplicial descent ones.
\end{obs}

Some direct consequences of the axioms are the following.

\begin{prop}\label{ExtraDeg}
The simple functor maps simplicial homotopy equivalences to weak equivalences. In particular, the following properties hold.\\
\textbf{\textit{i.}} If $f,g:X\rightarrow Y$ are simplicially homotopic maps, then $\mbf{s}(f)=\mbf{s}(g)$ in $\mc{D}[\mrm{E}^{-1}]$.\\
\textbf{\textit{ii.}} If $\alpha:X\rightarrow X_{-1}$ is an augmentation with an extra degeneracy then $\mbf{s}(\alpha)\in\mrm{E}$.
\end{prop}

\begin{proof} Let us see that $\mbf{s}(d_0^X:X\rightarrow X\boxtimes\Dl[1])$ is a weak equivalence for each $X$ in $\simp\mc{D}$. Note that $X\boxtimes\Dl[1]$ is the diagonal of $Z\in\simp\simp\mc{D}$ with $Z_{n,m}=\coprod_{\Dl[1]_n} X_m$. Therefore,
by (S2) $\mbf{s}\mrm{D}(Z)=\mbf{s}(X\boxtimes\Dl[1])$ and $\mbf{s}\mbf{s}(Z)$ are isomorphic in $\mc{D}[\mrm{E}^{-1}]$. By
(S1), there is a natural degreewise weak equivalence between $\mbf{s}(m\rightarrow Z_{\cdot,m})$ and $\mbf{s}(X)\boxtimes \Dl[1]$. So (S4) ensures that $\mbf{s}(X\boxtimes \Dl[1])$ and  $\mbf{s}(\mbf{s}(X)\boxtimes\Dl[1])$ are naturally isomorphic in $\mc{D}[\mrm{E}^{-1}]$.
Hence, by (S5), $\mbf{s}(d_0^X:X\rightarrow X\boxtimes \Dl[1])$ is in $\mrm{E}$. Since $d_0^X$ and $d_1^X:X\rightarrow X\boxtimes \Dl[1]$ have $s_0^X:X\boxtimes \Dl[1]\rightarrow X$ as common section, it follows that
$\mbf{s}(d_0^X)=\mbf{s}(d_1^X)$ in $\mc{D}[\mrm{E}^{-1}]$. This fact implies easily all the statements in the proposition.
 \end{proof}

Next result states that simplicial descent structures are inherited by diagram categories. The proof is straightforward, and is left to the reader.

\begin{prop}\label{DescensoFuntores} Let $I$ be a small category and let $(\mbf{s},\mu,\lambda)$ be a simplicial descent structure on $(\mc{D},\mrm{E})$. Then, the triple
$(\mbf{s}_{I\mc{D}},\mu_{I\mc{D}},\lambda_{I\mc{D}})$
defined objectwise is a simplicial descent structure on $(I\mc{D},I\mrm{E})$, where $I\mc{D}$ is the category of functors from $I$ to $\mc{D}$ and $I\mrm{E}=\{\alpha\mbox{ with }\alpha_i\in\mrm{E}\mbox{ for all }i\in I\}$. If $X:\simp\rightarrow I\mc{D}$ then
$(\mbf{s}_{I\mc{D}}(X))(i)=\mbf{s}(n\rightarrow X_n(i))$, and $(\mu_{I\mc{D}},\lambda_{I\mc{D}})$ is defined analogously.
\end{prop}

We will use of the notion of homotopical equivalence of relative categories of \cite[8.3 (ii)]{DHKS}.
In contrast to other homotopy theories based on the existence of cofibrations, such as
Quillen models, it holds that simplicial descent structures are closed by
homotopical equivalence of relative categories.

\begin{defi}\label{defiHomotEq} A \textit{homotopical equivalence} between the relative categories $(\mc{C},\mc{W})$ and $(\mc{D},\mrm{E})$ is
given by:\\
\textit{1.-} Functors $F:\mc{C}\rightarrow\mc{D}$ and $G:\mc{D}\rightarrow \mc{C}$ such that $F(\mc{W})\subset \mrm{E}$ and $G(\mrm{E})\subset\mc{W}$.\\
\textit{2.-} Zigzags of natural weak equivalences $\alpha : FG\dashrightarrow 1_{\mc{D}}$ and $\beta:GF\dashrightarrow 1_{\mc{C}}$.\\
We say that $(\mc{C},\mc{W})$ and $(\mc{D},\mrm{E})$ are \textit{homotopically equivalent} if there exists a homotopical equivalence between them.\\
This means that there exists an equivalence of categories between $\mc{C}[\mc{W}^{-1}]$ and $\mc{D}[\mrm{E}^{-1}]$ which is realized by weak equivalence-preserving functors
$F:\mc{C}\leftrightarrows\mc{D} :G$ and, in addition, the isomorphisms $FG\cong 1_{\mc{D}[\mrm{E}^{-1}]}$, $GF\cong 1_{\mc{C}[\mc{W}^{-1}]}$ are realized by zigzags of
natural weak equivalences.
\end{defi}

\begin{prop}\label{EquivHomot} Let $(\mc{C},\mc{W})$ and $(\mc{D},\mrm{E})$ be relative categories closed by finite coproducts, which are in addition homotopically equivalent.
If $(\mc{D},\mrm{E})$ is a simplicial descent category, then so is $(\mc{C},\mc{W})$.
\end{prop}

\begin{proof}
Let $(\mbf{s}^{\mc{D}},\mu^{\mc{D}},\lambda^{\mc{D}})$ be a simplicial descent structure on $(\mc{D},\mrm{E})$ and let $F:\mc{C}\leftrightarrows \mc{D}:G$,
$\alpha : FG\dashrightarrow 1_{\mc{D}}$,  $\beta:GF\dashrightarrow 1_{\mc{C}}$ be a homotopical equivalence between $(\mc{C},\mc{W})$ and $(\mc{D},\mrm{E})$.
We assume that $\beta_G \, G(\alpha^{-1}) = 1_G$ in $Fun(\mc{D},\mc{C})[\mc{W}^{-1}]$ and $F(\beta)\, \alpha_F^{-1} = 1_{F}$ in $Fun(\mc{C},\mc{D})[\mrm{E}^{-1}]$.
This is possible by lemma \ref{lemaEquivCatH}.
Let us see that $\mbf{s}^{\mc{C}} = G\mbf{s}^{\mc{D}}F : \simp\mc{C}\rightarrow \mc{C}$ is a simple functor for $(\mc{C},\mc{W})$.

(S1) We have that $\mbf{s}^{\mc{C}}$ preserves finite coproducts up to weak equivalence because $G$, $\mbf{s}^{\mc{D}}$ and $F$ do. We already know that $\mbf{s}^{\mc{D}}$ does.
On the other hand, $\mc{C}[\mc{W}^{-1}]$ and $\mc{D}[\mrm{E}^{-1}]$ are closed by finite coproducts, and they are preserved by $\gamma_{\mc{C}}:\mc{C}\rightarrow \mc{C}[\mc{W}^{-1}]$
and $\gamma_{\mc{D}}:\mc{D}\rightarrow \mc{D}[\mc{W}^{-1}]$. As $F:\mc{C}[\mc{W}^{-1}]\rightarrow \mc{D}[\mrm{E}^{-1}]$
is an equivalence of categories, it commutes with finite coproducts. This means that given $c,d\in\mc{C}$, the canonical morphism $\tau_F:F(c\sqcup d)\rightarrow F(c)\sqcup F(d)$ of
$\mc{D}$ becomes an isomorphism in $\mc{D}[\mrm{E}^{-1}]$. So $\tau_F\in\mrm{E}$ because $\mrm{E}$ is saturated. Arguing analogously for $G$, we get that $\mbf{s}^{\mc{C}}$ satisfies (S1).

(S2) Consider $Z\in\simp\simp\mc{C}$. Applying $F$ we obtain $F(Z) \in \simp\simp\mc{D}$, and we have the zigzag of natural weak equivalences
$\mu^{\mc{D}}_Z : \mbf{s}^{\mc{D}}F(\mrm{D}Z)\dashrightarrow \mbf{s}^{\mc{D}}\mbf{s}^{\mc{D}} F(Z)$. Then
$G(\mu^{\mc{D}}_Z):\mbf{s}^{\mc{C}}(\mrm{D}Z)\dashrightarrow G\mbf{s}^{\mc{D}}\mbf{s}^{\mc{D}} F(Z)$ is a zigzag of natural weak equivalences. On the other hand, applying $G\mbf{s}^{\mc{D}}$ to
$\alpha_{\mbf{s}^{\mc{D}} F(Z)}:FG\mbf{s}^{\mc{D}} F(Z)\dashrightarrow \mbf{s}^{\mc{D}} F(Z)$ we obtain the zigzag of natural weak equivalences $\rho_Z:G\mbf{s}^{\mc{D}}FG\mbf{s}^{\mc{D}} F(Z)=\mbf{s}^{\mc{C}}\mbf{s}^{\mc{C}}(Z)\dashrightarrow G\mbf{s}^{\mc{D}}\mbf{s}^{\mc{D}} F(Z)$.
Then, we set $\mu^{\mc{C}}_Z = \rho_Z^{-1} G(\mu^{\mc{D}}_Z)$.

(S3) Given $A\in\mc{M}$, we have $\lambda^{\mc{D}}_{F(A)} : \mbf{s}^{\mc{D}}c F(A) \rightarrow F(A)$. So $\lambda^{\mc{C}}_A:\mbf{s}^{\mc{C}}c(A)\rightarrow A$ is obtained from $G(\lambda^{\mc{D}}_{F(A)})$ and
$\beta_A:GF(A)\rightarrow A$. The compatibility between $\lambda^{\mc{C}}$ and $\mu^{\mc{C}}$ follows easily from the compatibility between $\lambda^{\mc{D}}$ and $\mu^{\mc{D}}$ and
from the equalities $\beta_G \, G(\alpha^{-1}) = 1_{G}$, $F(\beta)\, \alpha_F^{-1} = 1_{F}$.

(S4) $\mbf{s}^{\mc{C}}$ preserves degreewise weak equivalences because $\mbf{s}^{\mc{D}}$ does, $F(\mc{W})\subset\mrm{E}$ and $G(\mrm{E})\subset \mc{W}$.

(S5) Consider $A$ in $\mc{C}$. Since $F$ commutes with finite coproducts up to weak equivalence, for each simplicial finite set $K$ we have
a natural degreewise weak equivalence $F(A)\boxtimes K\rightarrow F (A\boxtimes K)$. It follows that
$\mbf{s}^{\mc{D}}( F(A)\boxtimes K )\rightarrow \mbf{s}^{\mc{D}} F (A\boxtimes K) $ is in $\mrm{E}$,
so $G\mbf{s}^{\mc{D}}( F(A)\boxtimes K )\rightarrow G\mbf{s}^{\mc{D}} F (A\boxtimes K)$ is in $\mc{W}$. Then $G\mbf{s}^{\mc{D}}( d_0^{F(A)})\sim \mbf{s}^{\mc{C}}(d_0^A)$, and
$\mbf{s}^{\mc{C}}(d_0^A)$ is in $\mc{W}$ since we know that $\mbf{s}^{\mc{D}}( d_0^{F(A)})$ is in $\mrm{E}$ and $G(\mrm{E})\subset \mc{W}$.
\end{proof}

Next lemma is a relative version of the fact that any pair $F:\mc{C}\leftrightarrows \mc{D}:G$ of inverse equivalences of categories
gives an `adjoint equivalence of categories' (see \cite[IV.4, Theorem 1]{ML}).

\begin{lema}\label{lemaEquivCatH}
Assume that $F:\mc{C}\leftrightarrows \mc{D}:G$, $\alpha : FG\dashrightarrow 1_{\mc{D}}$ and  $\beta:GF\dashrightarrow 1_{\mc{C}}$ form a homotopical equivalence between
$(\mc{C},\mc{W})$ and $(\mc{D},\mrm{E})$. Then there exists a zigzag $\beta':GF\dashrightarrow 1_{\mc{C}}$ of natural weak equivalences such that the compositions
$$F\stackrel{\alpha_F^{-1}}{\dashrightarrow} FGF \stackrel{F(\beta')}{\dashrightarrow} F \ \ \  \ \ G \stackrel{G(\alpha^{-1})}{\dashrightarrow} GFG \stackrel{\beta'_G}{\dashrightarrow} G$$
are equal to $1_{F}$ in $Fun(\mc{C},\mc{D})[\mrm{E}^{-1}]$ and to $1_G$ in $Fun(\mc{D},\mc{C})[\mc{W}^{-1}]$, respectively.
\end{lema}

\begin{proof} Since $\alpha$ is an isomorphism in $Fun(\mc{D},\mc{D})[\mrm{E}^{-1}]$ and $FG(\alpha)\, \alpha = \alpha_{FG} \, \alpha$ then  $FG(\alpha)=\alpha_{FG}$. Analogously, $GF(\beta)=\beta_{GF}$.
Note that $\beta':GF\dashrightarrow 1_{\mc{C}}$ in $Fun(\mc{C},\mc{C})[\mc{W}^{-1}]$ is determined by $\beta'_G : G F G\dashrightarrow G$. Indeed, since $\beta'$ is natural we must have the equality $\beta'\, GF (\beta) = \beta\, \beta'_{GF}$ in $Fun(\mc{C},\mc{C})[\mc{W}^{-1}]$. But $\beta$ and $GF (\beta)$ are isomorphisms in
$Fun(\mc{C},\mc{C})[\mc{W}^{-1}]$, so $\beta'$ is determined by $\beta'_G$. Therefore, we pick up the (unique)
$\beta':GF\dashrightarrow 1$ such that $\beta'_G = G(\alpha)$, that is $\beta' = \beta\, G(\alpha_{F})\, GF (\beta)^{-1}$. Then $\beta'$ is a zigzag of
natural weak equivalences. In addition $\beta'_G = \beta_G\, G(\alpha_{FG})\, GF (\beta_G)^{-1}=\beta_G\, GFG(\alpha)\, (\beta_{GFG})^{-1} = G(\alpha)$.\\
To finish it remains to see the equality $F(\beta') = \alpha_F$. Arguing as before, we have that $\gamma:F\dashrightarrow F$ in $Fun(\mc{C},\mc{D})[\mrm{E}^{-1}]$ is determined
by $\gamma_G$. Therefore it suffices to see that $F(\beta'_G) = \alpha_{FG}$.
But $\beta'_G = G(\alpha)$ implies $F(\beta'_G)=FG(\alpha)=\alpha_{FG}$.
\end{proof}

%%%%%%%%%%%%%%%%%%%%%%%%%%%%%%%%%%%%%%%%%%%%%%%%%%%%%%%%%%%%%%%%%%%%%%%%%%%%%%%%%%%%%%%%%%%%%%%%%%%%%%%%%%%%%%%%%%%%%%%%%%%%%%%%%%
%%%%%%%%%%%%%%%%%%%%%%%%%%%%%%%%%%%%%%%%%%%%%%%%%%%%%%%%%%%%%%%%%%%%%%%%%%%%%%%%%%%%%%%%%%%%%%%%%%%%%%%%%%%%%%%%%%%%%%%%%%%%%%%%%%
%%%%%%%%%%%%%%%%%                   CAPITULO
%%%%%%%%%%%%%%%%%%%%%%%%%%%%%%%%%%%%%%%%%%%%%%%%%%%%%%%%%%%%%%%%%%%%%%%%%%%%%%%%%%%%%%%%%%%%%%%%%%%%%%%%%%%%%%%%%%%%%%%%%%%%%%%%%%

\section{Examples}

In this section we exhibit some examples of simplicial and cosimplicial descent categories.
We begin with two general examples. In the first one, the weak equivalences form the biggest possible class: $\mrm{E}=\{$all morphisms$\}$. In the second one, $\mrm{E}$ is the smallest possible class: $\mrm{E}=\{$isomorphisms$\}$.

After that, we treat the classical examples of chain complexes and simplicial sets, which are useful to illustrate the axioms' meaning when the simple functor corresponds to the
`total complex' of a double chain complex on one hand, and the diagonal of a bisimplicial set on the other.

Further examples are deduced from theorem \ref{Modelos} in next section, where we prove that Quillen simplicial model categories are simplicial descent categories.
In the last part of this section we study the category of mixed Hodge complexes, as a remarkable example outside the setting of Quillen models. We see there that Deligne's
construction for cosimplicial mixed Hodge complexes provides a simple functor satisfying our axioms.

\begin{num}\label{EjGenerico} $\mbf{\mrm{E}}=\{$\textbf{all morphisms}$\}$.\\
Let $\mc{D}$ be a category with finite coproducts and initial object. Then, $(\mc{D},\mrm{E}=\{$all morphisms$\})$ admits a trivial simplicial descent structure $(\mbf{s},\mu,\lambda)$. The simple functor $\mbf{s}:\simp\mc{D}\rightarrow\mc{D}$ is $\mbf{s}(X)=X_k$ for a fixed $k\geq 0$, and $\lambda$ and $\mu$ are the identity natural transformations.
\end{num}

\begin{num} $\mbf{\mrm{E}}=\{$\textbf{isomorphisms}$\}$.\\
Let $\mc{C}$ be a category with finite colimits and initial object. Then, $(\mc{C},\mrm{E}=\{$isomorphisms$\})$ admits the following simplicial descent structure. The simple functor is $\mbf{s}=\mrm{eq} :\simp\mc{C}\rightarrow\mc{C}$, where $\mrm{eq}(X)$ denotes the coequalizer of
$$\xymatrix@C=15pt{X_0 && {\;} X_{1} \; \ar@<0.5ex>[ll]^-{d_0} \ar@<-0.5ex>[ll]_-{d_{1}} }$$
Note that $\mrm{eq}(X)$ agrees with the colimit of the whole diagram $X$. On the other hand, $\mu$ and $\lambda$ are the identity natural transformations.
\end{num}

\begin{num}\label{EjsSet} \textbf{Simplicial Sets and weak equivalences}.\\
Consider the class $\mc{W}$ of weak (homotopy) equivalences in $\simp Set$. Then, the diagonal $\mrm{D} : \simp\simp Set\rightarrow \simp Set$ endows $(\simp Set, \mc{W})$ with the simplicial descent structure $(\mrm{D}, \mu=id , \lambda=id)$.\\
Indeed, (S1), (S2) and (S3) are obvious. (S5) is a basic property of weak equivalences: they contain the homotopy equivalences. The remaining axiom (S4) is a well-known property:
\begin{quote}
(S4) Consider a map $F_{\cdot,\cdot}:X\rightarrow Y$ of bisimplicial sets such that for all $n\geq 0$, $F_{n,\cdot}$ is a weak equivalence.
Then $\mrm{D}(F)$, the diagonal of $F$, is again a weak equivalence.
\end{quote}
For a proof the reader may consult, for instance, \cite[proposition 1.9, p. 211]{GJ}. Also, it is possible to deduce the previous simplicial descent structure from theorem
\ref{Modelos}, or from proposition \ref{DlClosed},\textit{i}.
\end{num}

\begin{num}\label{EjCdc} \textbf{Positive chain complexes}.\\
Let $\mc{A}$ be an abelian category, and denote by $\cdcp{\mc{A}}$ the category of positive chain complexes. More concretely, $X\in\cdcp{\mc{A}}$ is a chain complex $\{X_n\}$ with $X_n=0$ for $n< 0$. Consider the class $\mrm{E}$ of weak equivalences  formed by the quasi-isomorphisms, that is, those chain maps inducing isomorphism on homology.\\
If $X=\{X_n , d_i , s_j \}$ is in $\simp \cdcp{\mc{A}}$, each $X_n$ is a
chain complex $\{ X_{n,p} , d_{X_n} \}_{p\in\mathbb{Z}}$. Hence
$X$ induces a double complex $K{X}$ with
$(K{X})_{n,p}=X_{n,p}$. The boundary maps are
$d_{X_n}:X_{n,p}\rightarrow X_{n,p-1}$ and
$\partial: X_{n,p}\rightarrow X_{n-1,p}$, $\partial=\sum_{i=0}^{n}(-1)^i d_i \,$.
The \textit{simple} functor $\mbf{s}:\simp \cdcp{\mc{A}}\rightarrow \cdcp{\mc{A}}$ is defined as $\mbf{s}(X)=\{$ total complex of $KX\}$, that is
$$\begin{array}{llr} (\mathbf{s}X)_q=\ds\bigoplus_{p+n=q} X_{n,p} & &
 d= \bigoplus (-1)^p\partial + d_{X_n}:\ds\bigoplus_{p+n=q} X_{n,p}\longrightarrow\ds\bigoplus_{p+n=q-1} X_{n,p}\end{array} $$
The following are well-known properties of $\mbf{s}$, which are inherited by those of the total complex functor.\\
(S1): $\mbf{s}$ is an additive functor.\\
(S2): Eilenberg-Zilber \cite[2.15]{DP}: If $Z\in\simp\simp \cdcp{\mc{A}}$, the
Alexander-Whitney map $\mu_Z:\mbf{s}\mrm{D}(Z)\rightarrow
\mbf{s}\mbf{s}(Z)$ and the `shuf{f}le' or Eilenberg-Zilber map
$\nu_Z:\mbf{s}\mbf{s}(Z)\rightarrow \mbf{s}\mrm{D}(Z)$ are inverse
homotopy equivalences.\\
In degree $n$, $(\mu_Z)_n$ is the sum of the maps
$Z(d^0\stackrel{j)}{\cdots}d^0,d^pd^{p-1}\cdots
d^{j+1}):Z_{p,p,q}\rightarrow Z_{i,j,q}$, $i+j=p$, $p+q=n$. The shuffle map is defined in degree $n$ by $(\nu_Z)_n = \bigoplus_{i+j=n}\nu_Z(i,j)$,  where
$$\nu_Z(i,j)=\sum_{(\alpha,\beta)}\epsilon(\alpha,\beta)Z(s^{\alpha_j}s^{\alpha_{j-1}}\cdots s^{\alpha_1},s^{\beta_i}s^{\beta_{i-1}}\cdots s^{\beta_1}):Z_{i,j}\rightarrow Z_{i+j,i+j}$$
The last sum is indexed over the $(i,j)$-shuffles $(\alpha,\beta)$, and $\epsilon(\alpha,\beta)$ is the sign of $(\alpha,\beta)$ \cite{EM}.\\
(S3): For each chain complex $A\in \cdcp{\mc{A}}$ there is a
natural splitting $\mbf{s}c(A)\cong A \oplus G$ where $G$ is contractible, and $\lambda_A$ is just the projection $A\oplus G\rightarrow A$.\\
(S4): If $X\in \simp\cdcp{\mc{A}}$ is such that $X_n$ is acyclic for all $n$, then
$\mbf{s}(X)$ is so. This is a well-known property of first-quadrant double complexes.
Property (S4) is easily deduced from this particular case.\\
(S5): Given a chain complex $A$, $\mbf{s}(A\boxtimes\Dl[1])$ is naturally homotopic to $cyl(A)$, the classical cylinder of $A$. See, for instance, \cite{W} for the precise definition of $cyl(A)$. I follows that if $e$ is a simplicial homotopy equivalence in $\simp\cdcp{\mc{A}}$ then $\mbf{s}(e)$ is a
homotopy equivalence in $\cdcp{\mc{A}}$, so (S5) holds.\\[0.1cm]
Therefore, for any abelian category $\mc{A}$, $(\cdcp{\mc{A}},\mrm{E})$ is a simplicial descent category with $(\mbf{s},\mu,\lambda)$ as defined above. Dually, the positive cochain complexes on $\mc{A}$, together with the quasi-isomorphisms as weak equivalences form a cosimplicial descent category, with the dual cosimplicial descent structure $(\mbf{s}:\Dl\cdcop{\mc{A}}\rightarrow \cdcop{\mc{A}},\mu,\lambda)$.
\end{num}
\begin{num}\label{EjHodge}{ \textbf{Mixed Hodge Complexes}}.\mbox{}\\
Next we define a category of mixed Hodge complexes and endow it with a structure
of cosimplicial descent category. The simple functor agrees with
Deligne's construction \cite[8.I.15]{DeIII} on objects.
It becomes a functor since the comparison morphisms of our mixed Hodge complexes are genuine filtered quasi-isomorphisms
instead of maps in the corresponding filtered derived category.\\
We denote by $\mathbb{Q}$ and $\mathbb{C}$ the respective categories
of $\mathbb{Q}$ and $\mathbb{C}$-vector spaces. Also, by $\mrm{CF}^+\mathbb{Q}$ and $\mrm{CF}^+\mathbb{C}$ we mean the respective categories of
filtered positive cochain complexes of $\mathbb{Q}$ and $\mathbb{C}$-vector spaces. All the filtrations are assumed to be biregular.
\end{num}

\begin{defi} A \textit{mixed Hodge complex} is the data
$((K_\mathbb{Q},\mrm{W}),(K_\mathbb{C},\mrm{W},\mrm{F}),\alpha)$,
where\\
\textit{1.-} $(K_\mathbb{Q},\mrm{W})\in \mrm{CF}^+\mathbb{Q}$ is such that $K_\mathbb{Q}$ has finite dimensional cohomology,
and $\mrm{W}$ is an increasing filtration.\\
\textit{2.-} $(K_\mathbb{C},\mrm{W},\mrm{F})$ is a positive bifiltered cochain complex of
$\mathbb{C}$-vector spaces, where $\mrm{W}$ (resp. $\mrm{F}$) is
an increasing (resp. decreasing) filtration, called the \textit{weight} (resp. \textit{Hodge}) filtration.\\
\textit{3.-} $\alpha$ is the data
$(\alpha_0,\alpha_1,(\widetilde{K},\widetilde{\mrm{W}}))$, where
$(\widetilde{K},\widetilde{\mrm{W}})$ is an object of
$\mrm{CF}^+\mathbb{C}$ and $\alpha_i$, $i=0,1$, is a filtered
quasi-isomorphism. That is, if $\mbf{Gr}_k:\mrm{CF}^+\mathbb{C}\rightarrow \cdcop{\mathbb{C}}$ denotes the graded functor then
$\mbf{Gr}_k(\alpha_i)$ is a quasi-isomorphism for all $k$ and $i=0,1$. Visually, $\alpha$ is
$$\xymatrix@M=4pt@H=4pt@C=35pt{ (K_\mathbb{C},\mrm{W}) & (\widetilde{K},\widetilde{\mrm{W}}) \ar[l]_-{\alpha_0} \ar[r]^-{\alpha_1} & (K_\mathbb{Q},\mrm{W})\otimes\mathbb{C} }$$
The following axiom must be satisfied\\
\textbf{(MHC)}
For each $n$ the boundary map of ${}_\mrm{W}\mbf{Gr}_nK_\mathbb{C}$ is compatible with the induced filtration $\mrm{F}$, and\\
$({}_\mrm{W}\mbf{Gr}_n H^k K_\mathbb{C},\mrm{F})$ is a Hodge structure of weight $n+k$. That is,
$$ {}_{\mrm{F}}\mbf{Gr}^p {}_{\overline{\mrm{F}}}\mbf{Gr}^q {}_\mrm{W}\mbf{Gr}_n H^k K_\mathbb{C} = 0\mbox{ for }p+q\neq n+k $$
\end{defi}

\begin{obs}\mbox{}\\
\textsc{I.} The filtered quasi-isomorphisms
have a calculus of fractions in the category of filtered
complexes up to filtered homotopy (see \cite[p.271]{I1}). Hence any
filtered quasi-isomorphism is represented by a zigzag as in \textit{3}.\\
\textsc{II.} Except for the $\mathbb{Z}$-part, a mixed Hodge complex as above is a mixed Hodge complex in the sense of Deligne, viewing $\alpha$ as an isomorphism in the filtered derived
category. Also, applying the decalage filtration to $\mrm{W}$ we get a mixed Hodge complex as defined in \cite{Hb} and \cite{B}.\\
\textsc{III.} We dropped the $\mathbb{Z}$-part of a mixed Hodge complex for simplicity, but all results in this section are also valid for mixed Hodge complexes with $\mathbb{Z}$-coef{f}icients.
\end{obs}

\begin{ej} \cite[8.I.8]{DeIII}
Let $j:U\rightarrow X$ be an open immersion  of complex smooth varieties,
where $X$ is proper and $Y=X\backslash U$ is a normal crossing divisor.\\
If $\mc{F}$ is a sheaf on $T$, set $R\Gamma (T,\mc{F}) = \Gamma(T, \mc{C}_{\mrm{God}}\mc{F})$, where $\mc{C}_{\mrm{God}}\mc{F}$ is the Godement resolution of $\mc{F}$. Analogously, if $\mc{F}$ is a bounded below complex of sheaves on $T$ (eventually filtered), set $R\Gamma (T,\mc{F}) = \Gamma(T, Tot(\mc{C}_{\mrm{God}}\mc{F}))$, where $Tot$ means the total complex of a double complex. The point is that $R\Gamma$ has values in the category of (filtered) complexes instead of the derived category, and the hypercohomology $H^\ast(T,\mc{F})$ may be computed as the cohomology of $R\Gamma (T,\mc{F})$.\\
Let $(\Omega_X\langle Y\rangle,\mrm{W},\mrm{F})$ be the logarithmic
De Rham complex of $X$ along $Y$ \cite[3.I]{DeII}. $\mrm{W}$ is the
so-called `weight filtration', and $\mrm{F}$ is the `Hodge
filtration', that is the filtration `b\^{e}te' associated with $\Omega_X\langle Y\rangle$.\\
Denote by $\mrm{W}$ the `canonical' filtration on $j_{\ast}\mathbb{Q}_U$, that is,
$\mrm{W}=\tau_{\leq} j_{\ast}\mathbb{Q}_U$.
A general argument shows that there is a zigzag of filtered quasi-isomorphisms connecting $R\Gamma (j_{\ast}\mathbb{Q}_U,\mrm{W})\otimes\mathbb{C}$ to $R\Gamma (\Omega_X\langle
Y\rangle,\mrm{W})$ (see \cite[p. 66]{H} or \cite[4.11]{PS}).
It is basically the result \cite[3.I.8]{DeII} connecting $\Omega_X\langle Y\rangle$ to $j_\ast \Omega_U$, together
with Poincar{\'e} lemma (that is, $\Omega_U$ is a resolution of the constant sheaf $\mathbb{C}_U$).
This zigzag may be reduced to a length 2 zigzag in a natural way (for instance, through the path object). Therefore $(R\Gamma (j_{\ast}\mathbb{Q},\mrm{W}),R\Gamma (\Omega_X\langle
Y\rangle,\mrm{W},\mrm{F}))$ is a mixed Hodge complex in the sense of previous definition.
\end{ej}

\begin{defi}
A morphism $(f_\mathbb{Q},f_\mathbb{C},\widetilde{f})\!
:\!((K_\mathbb{Q},\mrm{W}),(K_\mathbb{C},\mrm{W},\mrm{F}),\alpha)\rightarrow
((K'_\mathbb{Q},\mrm{W}'),(K'_\mathbb{C},\mrm{W}',\mrm{F}'),\alpha')$ of mixed Hodge complexes consists of morphisms
$f_\mathbb{Q}:(K_\mathbb{Q},\mrm{W})\rightarrow
(K'_\mathbb{Q},\mrm{W}')$ and
$f_\mathbb{C}:(K_\mathbb{C},\mrm{W},\mrm{F})\rightarrow
(K'_\mathbb{C},\mrm{W}',\mrm{F}')$ of (bi)filtered complexes.
If $\alpha$ and $\alpha'$ are the respective zigzags
$$\xymatrix@C=12pt{ (K_\mathbb{C},\mrm{W}) & (\widetilde{K},\widetilde{\mrm{W}}) \ar[l]_-{\alpha_0} \ar[r]^-{\alpha_1} & (K_\mathbb{Q},\mrm{W})\!\otimes\!\mathbb{C} &
 (K'_\mathbb{C},\mrm{W}') & (\widetilde{K}',\widetilde{\mrm{W}}') \ar[l]_-{\alpha'_0} \ar[r]^-{\alpha'_1} & (K'_\mathbb{Q},\mrm{W}')\!\otimes\!\mathbb{C}}$$
then $\widetilde{f}:(\widetilde{K},\widetilde{\mrm{W}})\rightarrow
(\widetilde{K}',\widetilde{\mrm{W}}')$ is a morphism of bifiltered complexes such that squares I and II in the diagram below
$$\xymatrix@H=4pt@C=23pt@R=17pt{(K_\mathbb{C},\mrm{W}) \ar[d]_{f_{\mathbb{C}}} & (\widetilde{K},\widetilde{\mrm{W}}) \ar[l]_-{\alpha_0} \ar[r]^-{\alpha_1}  \ar[d]_{\widetilde{f}}   & (K_\mathbb{Q},\mrm{W})\otimes\mathbb{C}  \ar[d]^{f_{\mathbb{Q}}\otimes\mathbb{C}} \\
                                (K'_\mathbb{C},\mrm{W}')\ar@{}[ru]|{\mrm{I}}   & (\widetilde{K}',\widetilde{\mrm{W}}') \ar[l]_-{\alpha'_0} \ar[r]^-{\alpha'_1} \ar@{}[ru]|{\mrm{II}} &
                                (K'_\mathbb{Q},\mrm{W}')\otimes\mathbb{C}}$$
commute.

In this way we obtain the category $\mc{H}dg$ of mixed Hodge complexes.
We consider the class of weak equivalences
${\mrm{E}}_{\mc{H}dg}=\{(f_\mathbb{Q},f_\mathbb{C},\widetilde{f})\; |\; f_\mathbb{Q}\mbox{ is a quasi-isomorphism }\mbox{in } \cdcop{\mathbb{Q}}\}$.
It follows from general Hodge theory that a weak equivalence induces an isomorphism between the corresponding mixed Hodge structures.
\end{defi}

Next we endow $(\mc{H}dg,\mrm{E}_{\mc{H}dg})$ with a cosimplicial descent structure.\\
\textit{Simple functor:}
If
${K}=((K_\mathbb{Q},\mrm{W}),(K_\mathbb{C},\mrm{W},\mrm{F}),\alpha)$
is a cosimplicial mixed Hodge complex, let ${\mbf{s}}_{\mc{H}dg}(K)$
be the mixed Hodge complex $\left((\mbf{s} (K_\mathbb{Q})
,\delta\mrm{W}),(\mbf{s}
(K_\mathbb{C}),\delta\mrm{W},\mbf{s}(\mrm{F})),\mbf{s}(\alpha)\right)$, where
$\mbf{s}$ denotes the simple of cochain complexes (see example \ref{EjCdc}) and
$\delta\mrm{W}$ is the diagonal filtration. More concretely
$$\begin{array}{cl}
 \mbf{s}(K_{\ast})^n=\ds\bigoplus_{p+q=n} K_{\ast}^{p,q} \ ;\ \ (\delta\mrm{W})_k(\mbf{s}(K_{\ast}))^n = \ds\bigoplus_{i+j=n}\mrm{W}_{k+i}K_{\ast}^{i,j} \ , & \mbox{ if }\ast\mbox{ is }\mathbb{Q}\mbox{ or }\mathbb{C}\\
 (\mbf{s}(\mrm{F}))^k (\mbf{s}(K_{\mathbb{C}}))^n = \ds\bigoplus_{p+q=n} \mrm{F}^k K_\mathbb{C}^{p,q} & {}
\end{array}$$
If $\alpha=(\alpha_0,\alpha_1,(\widetilde{K},\widetilde{\mrm{W}}))$
then $\mbf{s}(\alpha)$ denotes the zigzag
$$\xymatrix@M=4pt@H=4pt@C=34pt{
 (\mbf{s} (K_\mathbb{C}),\delta\mrm{W}) & (\mbf{s}(\widetilde{K}),\delta\widetilde{\mrm{W}}) \ar[l]_-{\mbf{s}(\alpha_0)} \ar[r]^-{\mbf{s}(\alpha_1)} & (\mbf{s} (K_\mathbb{Q}\otimes\mathbb{C}),\delta(\mrm{W}\otimes\mathbb{C})){\simeq }(\mbf{s} (K_\mathbb{Q}),\delta\mrm{W})\otimes \mathbb{C}}
$$
On morphisms, ${\mbf{s}}_{\mc{H}dg}(f_\mathbb{Q},f_\mathbb{C},\widetilde{f})=\left({\mbf{s}}(f_\mathbb{Q}),{\mbf{s}}(f_\mathbb{C}),{\mbf{s}}(\widetilde{f})\right)$.\\
\textit{Transformation}
$\mathbf{\lambda^{\mc{H}dg}}:1_{\mc{H}dg}\rightarrow
{\mbf{s}}_{\mc{H}dg}\comp c$ is
${\lambda}^{\mc{H}dg}_K=(\lambda^{\mathbb{Q}}_{K_{\mathbb{Q}}},\lambda^{\mathbb{C}}_{K_{\mathbb{C}}},\lambda^{\mathbb{C}}_{\widetilde{K}})$
induced by those $\lambda^\mathbb{Q}$ and $\lambda^{\mathbb{C}}$ of $\cdcop{\mathbb{Q}}$ and $\cdcop{\mathbb{C}}$ respectively.\\
\textit{Transformation} $\mathbf{\mu^{\mc{H}dg}}:\mbf{s}_{\mc{H}dg}\comp \mbf{s}_{\mc{H}dg}\rightarrow \mbf{s}_{\mc{H}dg}\comp \mrm{D}$
is defined analogously as $\mu^{\mc{H}dg}_K=(\mu^{\mathbb{Q}}_{K_{\mathbb{Q}}},\mu^{\mathbb{C}}_{K_{\mathbb{C}}},\mu^{\mathbb{C}}_{\widetilde{K}})$.

\begin{thm}\label{HodgeDescenso}
Deligne simple functor $\mbf{s}_{\mc{H}dg}:\Dl \mc{H}dg\rightarrow\mc{H}dg$ together with the transformations $\mu_{\mc{H}dg}$, $\lambda_{\mc{H}dg}$ defined above is a cosimplicial descent structure on $(\mc{H}dg,{\mrm{E}_{\mc{H}dg}})$.\end{thm}

\begin{proof} %
The proof is based on the fact that $\cdcop{\mathbb{Q}}$ is a cosimplicial descent category and the forgetful functor
$\mrm{U}:\mc{H}dg\rightarrow \cdcop{\mathbb{Q}}$ commutes with simple functors and the transformations $\lambda$ and $\mu$.\\
First of all, note that
${\mbf{s}}_{\mc{H}dg}=(\mbf{s},\delta,\mbf{s}):\Dl\mc{H}dg\rightarrow\mc{H}dg$
is indeed a functor.
Given ${K}\in\Dl\mc{H}dg$, then ${\mbf{s}}_{\mc{H}dg}(K)$ is a mixed
Hodge complex by \cite[8.I.15 i]{DeIII}, and ${\mbf{s}}_{\mc{H}dg}$
is functorial with respect to the morphisms of $\Dl\mc{H}dg$
by definition. Also, ${\mbf{s}}_{\mc{H}dg}$ is an additive functor, so (S1) holds.\\
To see (S2) and (S3), let $K=((K_\mathbb{Q},\mrm{W}),(K_\mathbb{C},\mrm{W},\mrm{F}),\alpha)$ be a mixed Hodge complex. Clearly $\lambda^\mathbb{Q}_{K_\mathbb{Q}}$,
$\lambda^\mathbb{C}_{K_\mathbb{C}}$ and $\lambda_{\widetilde{K}}^{\mathbb{C}}$ preserve the filtrations. Set $K=c(K)\in\Dl\mc{H}dg$. As the following diagram commutes in
$\mrm{CF}^+\mathbb{C}$
$$\xymatrix@H=4pt@C=19pt@R=17pt{(K_\mathbb{C},\mrm{W}) \ar[d]_{\lambda^{\mathbb{C}}_{K_\mathbb{C}}}   & (\widetilde{K},\widetilde{\mrm{W}}) \ar[l]_-{\alpha_0} \ar[r]^-{\alpha_1}  \ar[d]_{\lambda^{\mathbb{C}}_{\widetilde{K}}}                           & (K_\mathbb{Q},\mrm{W})\otimes\mathbb{C}  \ar[d]_{\lambda^{\mathbb{C}}_{K_\mathbb{Q}\otimes\mathbb{C}}} \ar[rd]^{\lambda^{\mathbb{Q}}_{K_\mathbb{Q}}\otimes\mathbb{C}}      &  \\
                                (\mbf{s}(K_\mathbb{C}),\delta(\mrm{W}))             & (\mbf{s}(\widetilde{K}),\delta(\widetilde{\mrm{W}})) \ar[l]_-{\mbf{s}(\alpha_0)} \ar[r]^-{\mbf{s}(\alpha_1)}  & (\mbf{s}(K_\mathbb{Q}\otimes\mathbb{C}),\delta(\mrm{W}\otimes\mathbb{C})) \ar[r]^-{\sim}   & (\mbf{s}(K_\mathbb{Q}),\delta(\mrm{W}))\otimes\mathbb{C}}$$
then
$\lambda^{\mc{H}dg}_K=(\lambda^{\mathbb{Q}}_{K_{\mathbb{Q}}},\lambda^{\mathbb{C}}_{K_{\mathbb{C}}},\lambda^{\mathbb{C}}_{\widetilde{K}})$
is a morphism in $\mc{H}dg$. Analogously
${\mu}^{\mc{H}dg}_K=(\mu^{\mathbb{Q}}_{K_{\mathbb{Q}}},\mu^{\mathbb{C}}_{K_{\mathbb{C}}},\mu^{\mathbb{C}}_{\widetilde{K}})$ is a morphism in $\mc{H}dg$. Since
$\lambda^{\mathbb{Q}}_{K_{\mathbb{Q}}}$ and $\mu^{\mathbb{Q}}_{K_{\mathbb{Q}}}$ are quasi-isomorphisms in $\cdcop{\mathbb{Q}}$, then ${\lambda}^{\mc{H}dg}_K$ and
${\mu}^{\mc{H}dg}_K$ are in $\mrm{E}_{\mc{H}dg}$.\\
Axiom (S4) is clear since ${\mbf{s}}_{\mc{H}dg}(f_\mathbb{Q},f_\mathbb{C},\widetilde{f})=\left({\mbf{s}}(f_\mathbb{Q}),{\mbf{s}}(f_\mathbb{C}),{\mbf{s}}(\widetilde{f})\right)$ and (S4)
holds in $\cdcop{\mathbb{Q}}$. Finally, given a simplicial finite set $L$, and a mixed Hodge complex $K$ then $\mrm{U}(K\boxtimes L)=\mrm{U}(K)\boxtimes L$, so (S5) holds
since it holds for cochain complexes.
\end{proof}

%%%%%%%%%%%%%%%%%%%%%%%%%%%%%%%%%%%%%%%%%%%%%%%%%%%%%%%%%%%%%%%%%%%%%%%%%%%%%%%%%%%%%%%%%%%%%%%%%%%%%%%%%%%%%%%%%%%%%%%%%%%%%%%%%%
%%%%%%%%%%%%%%%%%%%%%%%%%%%%%%%%%%%%%%%%%%%%%%%%%%%%%%%%%%%%%%%%%%%%%%%%%%%%%%%%%%%%%%%%%%%%%%%%%%%%%%%%%%%%%%%%%%%%%%%%%%%%%%%%%%
%%%%%%%%%%%%%%%%%                   CAPITULO
%%%%%%%%%%%%%%%%%%%%%%%%%%%%%%%%%%%%%%%%%%%%%%%%%%%%%%%%%%%%%%%%%%%%%%%%%%%%%%%%%%%%%%%%%%%%%%%%%%%%%%%%%%%%%%%%%%%%%%%%%%%%%%%%%%

\section{Homotopy colimits in simplicial model categories}

In this section we prove that the Bousfield-Kan formula for the homotopy colimit in a simplicial model category $(\mc{M},\mc{W})$
gives rise to a simplicial descent structure on the subcategory of cofibrant objects of $\mc{M}$.
In case $\mc{M}$ has functorial cofibrant replacements and $\mc{W}$ is closed by finite coproducts, the `corrected' Bousfield-Kan homotopy colimit
does induce a simplicial descent structure on all $\mc{M}$. The dual results hold for Bousfield-Kan homotopy limits and cosimplicial descent
structures.

We restrict ourselves to simplicial model categories because in this case homotopy limits and colimits are easier to define and to deal with.
In a model category which is not necessarily simplicial one defines Bousfield-Kan homotopy limits (resp. colimits) through the choice of a simplicial (resp. cosimplicial)
frame (see \cite{H}). The results given here also work in this general setting, although the proofs become more technical.

In the cubical case, the connection between simplicial model categories and (cubical) homological descent categories, which are developed in \cite{GN},
is studied in \cite{Ru}.

We assume the reader is familiar with simplicial model categories and Bousfield-Kan homotopy limits and colimits. We refer the reader to \cite{H} for definitions and proofs.\\

Given a functor $F:J\rightarrow I$ and an object $x\in I$, recall that the \textit{overcategory} $F\downarrow x$ has as objects the maps $F(y)\rightarrow x$ of $I$. We denote by $I\downarrow x$ the overcategory $1_{I}\downarrow x$. $F$ is called \textit{homotopy left cofinal} if the simplicial set $\mrm{B}(F\downarrow x)$  is contractible for each $x\in I$. Here $\mrm{B}(\mc{C})$ denotes the nerve of $\mc{C}$. $F$ is \textit{homotopy right cofinal} if $F^{op}:{J}^{op}\rightarrow I^{op}$ is homotopy left cofinal.

\begin{defi}
Let $X:I\rightarrow\mc{M}$ be a functor from a small category $I$ to a simplicial model category $(\mc{M},\otimes)$.
Consider the bifunctor $X\otimes \mrm{B}(I^{op}\downarrow \cdot):I\times I^{op}\rightarrow \mc{M}$; $(c,d)\mapsto X(c)\otimes \mrm{B}(I^{op}\downarrow d)$.
Then, the \textit{homotopy colimit} of $X$, $\mrm{hocolim}_{I} X$
(or $\mrm{hocolim}\, X$ if $I$ is understood), is the coend \cite[IX.6]{ML}
$$\mrm{hocolim}_{I} X = \int^c X(c)\otimes \mrm{B}(I^{op}\downarrow c)$$
A functor $F:J\rightarrow I$ induces a natural map $\mrm{hocolim}_{J}F^\ast X\rightarrow \mrm{hocolim}_{I}X$, defined by the maps $1\otimes \mrm{B}(F): X(F(b))\otimes  \mrm{B}(J^{op}\downarrow b)\rightarrow X(F(b))\otimes  \mrm{B}(I^{op}\downarrow F(b))$.
\end{defi}

We will use the following property of $\mrm{hocolim}$
\begin{quote}
\cite[19.6.13]{H} If $F:J\rightarrow I$ is homotopy right cofinal and $X:I\rightarrow\mc{M}$, then the induced map $\mrm{hocolim}_{J} F^\ast X\rightarrow \mrm{hocolim}_{I}X$ is a weak equivalence.
\end{quote}

\begin{thm}\label{Modelos} Let $\mc{M}$ be a simplicial model category with weak equivalences $\mc{W}$. Denote by $\mc{M}_c$ and $\mc{M}_f$ the respective
subcategories of cofibrant and fibrant objects of $\mc{M}$. Then, $(\mc{M}_c,\mc{W},\mrm{hocolim})$ is a simplicial descent category.
Dually, $(\mc{M}_f,\mc{W},\mrm{holim})$ is a cosimplicial descent category.
\end{thm}

\begin{proof} It is well known that cofibrant objects are closed by coproducts, and that the coproduct of two weak equivalences between cofibrant objects is again a weak equivalence. On the other
hand, $\mc{W}$ is saturated in $\mc{M}$ by \cite[8.3.9]{H}, so it is saturated in $\mc{M}_c$ as well.
To see (S1), we have that if $X:\simp\rightarrow\mc{M}$ is objectwise cofibrant then $\mrm{hocolim}_{\simp} X$ is cofibrant by \cite[18.5.2]{H}.
It follows that $\mrm{hocolim}$ is a functor $\simp \mc{M}_c\rightarrow \mc{M}_c$, and it is clear from the definition that it preserves finite coproducts.\\
The transformations $\lambda$ and $\mu$ are easily defined using that $\mrm{hocolim}_{I}$ is natural in $I$. Indeed, $\lambda$ is obtained from $l:\simp\rightarrow \ast$. It induces the
map $\lambda_M:\mrm{hocolim}_{\simp}c(M)\rightarrow M=\mrm{hocolim}_{\ast}M$ for each $M\in\mc{M}_c$. We have $\lambda_M\in\mc{W}$ because $l$ is homotopy right cofinal ($\simp$
has an initial object, and therefore contractible nerve).\\
On the other hand, the diagonal $d:\simp\rightarrow \simp\times\simp$ induces $\mrm{hocolim}_{\simp}\mrm{D}(Z)\rightarrow\mrm{hocolim}_{\simp\times\simp}Z$ for each
$ Z\in \simp\simp\mc{M}_c$. It gives us a transformation $\mu_Z:\mrm{hocolim}_{\simp}\mrm{D}(Z)\rightarrow\mrm{hocolim}_{\simp}\mrm{hocolim}_{\simp}Z$, since the Fubini property of
$\mrm{hocolim}$ ensures that $\mrm{hocolim}_{\simp\times\simp}$ and $\mrm{hocolim}_{\simp}$ applied twice are the same. But $d$ is homotopy right cofinal (see \cite[lemma 5.33]{T}), so
$\mu_Z\in\mc{W}$. Now, the compatibility between $\lambda$ and $\mu$ holds trivially since the composition of the diagonal $d:\simp\rightarrow\simp\times\simp$ with
$l\times 1:\simp\times \simp\rightarrow\simp$ (resp. with $1\times l$) is the identity.
So (S2) and (S3) are satisfied. (S4) is the property of $\mrm{hocolim}$ known as `homotopy invariance' (see \cite[18.5.3]{H}): if $X,Y:I\rightarrow\mc{M}$ are objectwise cofibrant and
$\tau:X\rightarrow Y$ is objectwise a weak equivalence, then  $\mrm{hocolim}_{I} \tau\in\mc{W}$.

To finish, it remains to prove (S5). If $A\in\mc{M}_c$, let us see that $\mrm{hocolim} (d_0^A:A\rightarrow A\boxtimes \Dl[1])$ is in $\mc{W}$. The key point is that there is a natural
weak equivalence $\mrm{hocolim} (A\boxtimes\Dl[1])\rightarrow A\otimes \Dl[1]$, where $\otimes$ is the internal action given by the simplicial model structure. To see this, note that
since $A$ is cofibrant, then for any simplicial set $K$, $A\boxtimes K$ is Reedy cofibrant. In this case the Bousfield-Kan map gives a weak equivalence between
$\mrm{hocolim} (A\boxtimes K)$ and $|A\boxtimes K|$, the geometric realization of $A\boxtimes K$ \cite[18.7.4]{H}. But  by \cite[(3.7), p. 385]{GJ} (or by direct computation), it holds
that $|A\boxtimes K|\simeq A\otimes K$. Finally, since $d^0:\Dl[0]\rightarrow\Dl[1]$ is a weak equivalence in $\simp Set$, then $A\otimes d^0:A\rightarrow A\otimes \Dl[1]$ is a weak
equivalence in $\mc{M}$ (see \cite[9.3.9]{H}). We conclude that $\mrm{hocolim} (d_0^A) = \mrm{hocolim}(A\boxtimes d^0)\sim A\otimes d^0$ is also a weak equivalence.
\end{proof}

\begin{obs}\label{NoQ} The converse of theorem \ref{Modelos} does not hold. Indeed, there are examples of simplicial descent categories whose homotopy category can not be equivalent to
the homotopy category of any Quillen model category.

For instance, there are abelian categories $\mc{A}$ whose associated (positive) derived category $\mrm{D}_+(\mc{A})$ does not have small hom's. There is an explicit example of
such $\mc{A}$ due to P. Freyd. It consists of the abelian category of small $R$-modules where $R$ is the `big' ring of polynomials on a proper class of variables, and with coefficients
in $\mathbb{Z}$.
 As we have seen in example \ref{EjCdc}, $(\mc{D}=\cdcp{\mc{A}}, \mrm{E}=\{$quasi-isomorphisms$\})$ is a simplicial descent category, and $\mc{D}[\mrm{E}^{-1}]=\mrm{D}_+(\mc{A})$.
But $\mathbb{Z}$ endowed with the trivial $R$-module structure is an object of $\mc{A}$ with a proper class of submodules. See \cite{F}, or \cite{CN}, for more details. It turns out that
$\mrm{Ext}^1(\mathbb{Z},\mathbb{Z})\simeq \mrm{Hom}_{\mrm{D}_+(\mc{A})}(\mathbb{Z},\Sigma\mathbb{Z})$ is a proper class. Hence, $\mc{D}[\mrm{E}^{-1}]$ can not be equivalent to the homotopy
category $\mc{M}[\mc{W}^{-1}]$ of any Quillen model category $(\mc{M},\mc{W})$, since in this case the morphisms in $\mc{M}[\mc{W}^{-1}]$ between any two fixed objects would form a
small set.
\end{obs}

The previous theorem provides a wide class of examples of simplicial and cosimplicial descent categories. For instance, we recover example \ref{EjsSet}.

To finish this section, we show that it is also possible to induce simplicial descent structures on  all $\mc{M}$, not only on the cofibrant objects. The price
is to `correct' the Bousfield-Kan homotopy colimit by composing it with an objectwise cofibrant replacement.

\begin{cor}\label{ModelosCompl} Let $(\mc{M},\mc{W})$ be a simplicial model category. Assume that $\mc{M}$ admits functorial factorizations, and choose functorial replacements
$Q_c : \mc{M}\rightarrow \mc{M}_c$, $Q_f : \mc{M}\rightarrow \mc{M}_f$.\\
\textbf{\textit{i.}} If $\mc{W}$ is closed by finite coproducts, then $(\mc{M},\mc{W})$ is a simplicial descent category with
simple functor $\mbf{s}=\mrm{hocolim}\, Q_c : \simp\mc{M}\rightarrow \mc{M}$.\\
\textbf{\textit{ii.}} If $\mc{W}$ is closed by finite products, then $(\mc{M},\mc{W})$ is a cosimplicial descent category with
simple functor $\mbf{s}=\mrm{holim}\, Q_f : \Dl\mc{M}\rightarrow \mc{M}$.
\end{cor}

\begin{proof}
By duality, it suffices to see the first part. Since a cofibrant functorial replacement $Q_c:\mc{M}\rightarrow \mc{M}_c$ and the inclusion
$i:\mc{M}_c\rightarrow \mc{M}$ form a homotopical equivalence between $(\mc{M},\mc{W})$ and $(\mc{M}_c,\mc{W})$, the result is a consequence
of previous theorem and proposition \ref{EquivHomot}.
\end{proof}

%%%%%%%%%%%%%%%%%%%%%%%%%%%%%%%%%%%%%%%%%%%%%%%%%%%%%%%%%%%%%%%%%%%%%%%%%%%%%%%%%%%%%%%%%%%%%%%%%%%%%%%%%%%%%%%%%%%%%%%%%%%%%%%%%%
%%%%%%%%%%%%%%%%%%%%%%%%%%%%%%%%%%%%%%%%%%%%%%%%%%%%%%%%%%%%%%%%%%%%%%%%%%%%%%%%%%%%%%%%%%%%%%%%%%%%%%%%%%%%%%%%%%%%%%%%%%%%%%%%%%

\section{$\Dl$-closed classes and Brown categories of cofibrant objects}\label{appendix}

We begin by reminding the reader the definition of $\Dl$-closed class, which is developed in \cite{Vo}.

\begin{defi}\label{defiCC}
Let $\mc{C}$ be a category. A class $\mc{W}$ of morphisms in $\simp\mc{C}$ is called $\Dl$\textit{-closed} if it satisfies the following three properties\\[0.1cm]
\textit{1.-} The class $\mc{W}$ contains the simplicial homotopy equivalences.\\[0.1cm]
\textit{2.-} $\mc{W}$ satisfies the 2-out-of-3 property.\\[0.1cm]
\textit{3.-} If $F=F_{\cdot,\cdot}:Z_{\cdot,\cdot}\rightarrow T_{\cdot,\cdot}$ is a map of bisimplicial objects in $\mc{C}$ such that
$F_{n,\cdot}\in\mc{W}$ (or $F_{\cdot,n}\in\mc{W}$) for all $n\geq 0$, then the diagonal of $F$, $\mrm{D}(F)$, is in $\mc{W}$.\\[0.1cm]
Assume moreover that  $\mc{C}$ has finite coproducts. Then $\mc{W}$ is $(\Dl,\amalg_{< \infty})$\textit{-closed} if it is $\Dl$-closed, and\\[0.1cm]
\textit{4.-} If $F,G \in\mc{W}$, then $F\sqcup G\in\mc{W}$.
\end{defi}

\begin{prop}\label{DlClosed}\mbox{}\\
\vspace{-0.5cm}
\begin{compactitem}
\item[\textbf{\textit{i}.}] Consider a category $\mc{C}$ with finite coproducts, and a saturated class $\mc{W}$ of morphisms in $\simp\mc{C}$. Then, the following are equivalent.
\begin{compactitem}
\item[1.] $\mc{W}$ is $(\Dl,\amalg_{< \infty})$-closed.
\item[2.] $(\simp\mc{C},\mc{W},\mrm{D}:\simp\simp\mc{C}\rightarrow\simp\mc{C})$ is a simplicial descent category.
\end{compactitem}
\item[\textbf{\textit{ii}.}] Given a simplicial descent category $(\mc{D},\mrm{E},\mbf{s})$, then
$\mc{S}=\mbf{s}^{-1}\mrm{E}$ is $(\Dl,\amalg_{< \infty})$-closed.
\end{compactitem}
\end{prop}

\begin{obs}
The saturated class $\mc{S}=\mbf{s}^{-1}\mrm{E}$ does not depend on the simple functor $\mbf{s}$. If $\mbf{s'}$ is another simple functor making $(\mc{D},\mrm{E})$ a simplicial descent category,
then $\mbf{s}^{-1}\mrm{E}=\mbf{s'}^{-1}\mrm{E}$. This fact follows from corollary \ref{unicSimple}.
\end{obs}

\noindent Before giving the proof of previous proposition, let us introduce some notations.\\
\hspace{0.2cm}\textmd{-} by $\simp\mrm{E}$ we mean the class of maps $f:X\rightarrow Y$ in $\simp\mc{D}$ such that $f_n\in\mrm{E}$ for all $n\geq 0$.\\
\hspace{0.2cm}\textmd{-} $\mbf{s}:\simp\mc{D}\rightarrow \mc{D}$ induces $\simp\mbf{s},\,\mbf{s}\simp :\simp\simp\mc{D}\rightarrow \simp\mc{D}$, defined as follows. If $Z$ is in $\simp\simp\mc{D}$,
    $$\begin{array}{l}
    (\simp\mbf{s})(Z)_n = \mbf{s}(m\rightarrow Z_{n,m})\\
    (\mbf{s}\simp)(Z)_n = \mbf{s}(m\rightarrow Z_{m,n})
    \end{array}$$

\begin{lema}\label{Wclass}\mbox{}\\
\textbf{\textit{i}.} $\simp\mrm{E} \subset \mc{S}$ in $\simp\mc{D}$.\\
\textbf{\textit{ii}.} $(\simp\mbf{s})^{-1}\mc{S}=(\mbf{s}\simp)^{-1}\mc{S}=\mrm{D}^{-1}\mc{S}$ in $\simp\simp\mc{D}$.
\end{lema}

\begin{proof} Part \textit{i} is just the exactness axiom (S4). Part \textit{ii} is an easy consequence of (S2). Indeed, let $F:Z\rightarrow T$ be a map in $\simp\simp\mc{D}$. Axiom (S2)
 produces a natural isomorphism $\mbf{s}\mrm{D}(F) \simeq \mbf{s}(\simp\mbf{s})(F)$ in $\mc{D}[\mrm{E}^{-1}]$. Then $\mrm{D}(F)\in\mbf{s}^{-1}\mrm{E}=\mc{S}$ if and only if
$(\simp\mbf{s})(F)$ is. On the other hand, if $F'$ is the bisimplicial map given by $F'_{n,m}=F_{m,n}$, then $\mbf{s}\mrm{D}(F) = \mbf{s}\mrm{D}(F')\simeq \mbf{s}(\simp\mbf{s})(F') =
\mbf{s}(\mbf{s}\simp)(F)$. Therefore, $\mrm{D}(F)\in\mc{S}$ if and only if $(\mbf{s}\simp)(F)$ is.
\end{proof}

\begin{proof}[Proof of proposition \ref{DlClosed}]
Let $\mc{C}$ and $\mc{W}$ be as in \textit{i}, and assume that \textit{i.1} holds. Then $(\mrm{D},\mu=id, \lambda=id)$ is a simplicial descent structure on
$(\simp\mc{C},\mc{W})$. Indeed, in this case the only non-trivial axioms are (S4) and (S5), which correspond to properties \textit{3} and \textit{1} of $(\Dl,\amalg_{< \infty})$-closed class, respectively.\\
Conversely, assume now that $(\simp\mc{C},\mc{W},\mrm{D})$ is a simplicial descent category.
Since $\mc{W}$ is saturated and closed by finite coproducts, then properties \textit{2} and \textit{4} of $(\Dl,\amalg_{< \infty})$-closed class hold. Let us see \textit{1}. Given $K\in\simp\mc{C}$, then $\mrm{D}(K\boxtimes\Dl[1])$ is by definition
equal to the usual cylinder object of $K$, $K\times \Dl[1]$. Then, property \textit{1} follows from (S5) and the 2-of-3 property of $\mc{W}$. To finish, \textit{3} is just the exactness axiom for $\mrm{D}$.
Indeed, given $F\in\simp\simp\mc{C}$ with $F_{n,\cdot}\in\mc{W}$ for all $n$, then $\mrm{D}(F)\in\mc{W}$ by (S4). If $F_{\cdot,n}\in\mc{W}$ for all $n$, let $F'\in\simp\simp\mc{C}$ be $F'_{n,m}=F_{m,n}$.
Then $\mrm{D}(F)=\mrm{D}(F')\in\mc{W}$.\\
Assume now that $(\mc{D},\mrm{E})$ is a simplicial descent category, and $\mbf{s}$ is a simple functor. Property \textit{1} of $(\Dl,\amalg_{< \infty})$-class for $\mc{S}$ is just proposition \ref{ExtraDeg}. On the other hand, $\mc{S}$ is saturated since $\mrm{E}$ is, so property \textit{2}
holds. To see \textit{3}, consider a bisimplicial map $F$ such that $F_{n,\cdot}\in\mc{S}$ (the case $F_{n,\cdot}\in\mc{S}$ is analogous). This means that $\mbf{s}(n\rightarrow F_{n,m})\in\mrm{E}$ for all $m\geq 0$,
or equivalently, that $(\mbf{s}\simp)(F)\in\simp\mrm{E}$. But $\simp\mrm{E}\subset \mc{S}$, so $F\in (\mbf{s}\simp)^{-1}\mc{S}=\mrm{D}^{-1}\mc{S}$ by lemma \ref{Wclass}. Hence $\mrm{D}(F)\in\mc{S}$. Finally,
\textit{4} follows directly from (S1) and the fact that $\mrm{E}$ is closed by finite coproducts.
\end{proof}

We describe in the remaining part of this section a
Brown structure of cofibrant objects on $(\simp\mc{C},\mc{W})$, where $\mc{W}$ is a $(\Dl,\amalg_{< \infty})$-closed class. We assume the reader is familiar with sections
1, 2 and 4 of \cite{Br}.

\begin{defi}\cite{Vo}
A map $F:X\rightarrow Y$ in $\simp\mc{C}$ is a \textit{termwise coprojection} if for each $n\geq 0$, there exists $A^{(n)}\in\mc{C}$
such that
$$\xymatrix@M=4pt@H=4pt@C=25pt@R=10pt{ X_n \ar[r]^{F_n} \ar[rd] &  Y_n \ar@{-}[d]_{\wr} \\
				                         &  X_n \sqcup A^{(n)} }$$
where $X_n\rightarrow X_n\sqcup A^{(n)}$ is the canonical map.\\
Note that, in this case, given any other map $G:X\rightarrow Z$, the pushout of $F$ and $G$ always exists in $\simp\mc{C}$. We denote it by $Y\cup_X Z$.\\[0.1cm]
Given a diagram $(Q)$: $Z\stackrel{G}{\leftarrow}X\stackrel{F}{\rightarrow} Y$ in $\simp\mc{C}$, define $K(F,G)$ by the pushout
$$\xymatrix@M=4pt@H=4pt@R=12pt@C=30pt{
 X\sqcup X \ar[d]_{F\sqcup G}\ar[r]^-{(d_0^X, d_1^X)} &  X\boxtimes \Dl[1] \ar[d]_-{}\\
 Y\sqcup Z\ar[r] & K(F,G)}$$
If $G=1_X : X\rightarrow X$, then $K(F,G)$ is called the \textit{simplicial cylinder} of $F$, and denoted by $\mathtt{Cyl}(F)$.
\end{defi}

Note that $ K(F,G)$ and $Y\sqcup Z\rightarrow  K(F,G)$ are natural on $(F,G)$.

\begin{lema}\label{Pushouts}\mbox{}\\
\textbf{\textit{i}.} If $F$ is a termwise coprojection, then the natural map $K(F,G)\rightarrow Y \cup_X Z $ is in $\mc{W}$.\\
\textbf{\textit{ii}.} Given a map of diagrams $(f_Z,f_X,f_Y):(Q)\rightarrow (Q')$ such that $f_Z,f_X$ and $f_Y$ are in $\mc{W}$, then the induced map
$K(F,G)\rightarrow K(F',G')$ is in $\mc{W}$ as well.\\
\textbf{\textit{iii}.} The natural maps $Y\rightarrow \mathtt{Cyl}(F)\rightarrow Y$ are inverse simplicial homotopy equivalences.
\end{lema}

\begin{proof}
The result is a consequence of lemmas 2.10, 2.11 and 2.9 in \cite{Vo}.
\end{proof}

We will need a notion of distinguished triangles in $\simp\mc{C}$ defined through the cone functor rather than through
cofibrations. The reason is that we want to transfer these constructions from $\simp\mc{C}$ to $\mc{C}$, and we do not have a notion of cofibrations in $\mc{C}$,
while we do have an induced cone functor.

\begin{defi}\label{DefiCono} Assume moreover that $\mc{C}$, and hence $\simp\mc{C}$, is pointed with zero object $\ast$. Consider a map $F:X\rightarrow Y$ in $\simp\mc{C}$, and denote by $\ast_X:X\rightarrow \ast$ the trivial map.\\
The \textit{simplicial cone} of $F$, denoted $\mathtt{Cone}(F)$, is by definition equal to $K(F,{\ast_X})$.\\
The \textit{simplicial suspension} of $X$, denoted $\Lambda (X)$, is the simplicial cone of $\ast_X$, that is $\Lambda(X)=K(\ast_X,\ast_X)$.\\
The \textit{cone sequence} induced by $F$  is the sequence in $\simp\mc{D}$ given by
\begin{equation}\label{SucCono} X\stackrel{F}{\rightarrow} Y \rightarrow \mathtt{Cone}(F) \end{equation}
\end{defi}

\noindent We define the natural map
\begin{equation}\label{Coaction}\mathtt{a}_F:\mathtt{Cone}(F)\rightarrow \mathtt{Cone}(F)\sqcup \Lambda(X)\end{equation}
in $(\simp\mc{C})[\mc{W}^{-1}]$ as follows. Set $\Theta(F)=\mathtt{Cone}\left(X {\rightarrow}\, \mathtt{Cyl}(F)\right)$.
The commutative square
$$\xymatrix@M=4pt@H=4pt@C=25pt@R=10pt{X\ar[d]\ar[r] & \mathtt{Cyl}(F) \ar[d]\\
                                      X  \ar[r] &  Y}$$
induces $\alpha: \Theta(F)\rightarrow\mathtt{Cone}(F)$, which is in $\mc{W}$.
On the other hand, using $\mathtt{Cyl}(F)\rightarrow \mathtt{Cone}(F)$ we get $\beta: \Theta(F)\rightarrow \mathtt{Cone}\left( X\rightarrow \mathtt{Cone}(F)\right)$, where
$X\rightarrow \mathtt{Cone}(F)$ is the trivial map factoring through $\ast$. Therefore $\mathtt{Cone}\left( X\rightarrow \mathtt{Cone}(F)\right)$ is canonically isomorphic to
$\mathtt{Cone}(F)\sqcup \Lambda(X)$. We define $\mathtt{a}_F := \beta\circ \alpha^{-1}$.\\

\noindent If $Y=\ast$, we deduce
\begin{equation}\label{CogrObj}\mathtt{a}_X:\Lambda(X)\rightarrow \Lambda(X)\sqcup \Lambda(X)\end{equation}
which endows $\Lambda(X)$ with a cogroup structure, as we will see below. In addition, for a general map $F$, $\mathtt{a}_F$ gives a coaction of $\Lambda(X)$ on $\mathtt{Cone}(F)$.

\begin{defi}
A \textit{simplicial distinguished triangle} in $(\simp\mc{C})[\mc{W}^{-1}]$ is a pair
$$\left( X\rightarrow Y \rightarrow Z  \, , \, Z\rightarrow Z\sqcup \Lambda(X) \right)$$
isomorphic in $(\simp\mc{C})[\mc{W}^{-1}]$ to the one induced by some $F:X\rightarrow Y$ in $\simp\mc{C}$ as in (\ref{SucCono}) and (\ref{Coaction}).
\end{defi}

Although the following result is not stated explicitly in \cite{Vo}, its proof is contained there.

\begin{prop}\label{VovWald}\mbox{}\\
\textbf{\textit{i.}} Let $\mc{C}$ be a category with finite coproducts and $\mc{W}$ a $(\Dl,\amalg_{< \infty})$-closed class of $\simp\mc{C}$. Then
$(\mc{C},\mc{W},Cof)$ is a Brown category of cofibrant objects, where $Cof=\{$termwise coprojections$\}$.\\
\textbf{\textit{ii.}} If in addition $\mc{C}$ is pointed, then a pair $\left( X\rightarrow Y \rightarrow Z  \, , \, Z\rightarrow Z\sqcup \Lambda(X) \right)$ is a simplicial distinguished triangle if and only if it is a cofibration sequence in the sense of \emph{\cite{Br}}.
\end{prop}

\begin{proof}
To see the first part, the only non-trivial properties to check are the pushout and cylinder axioms. The pushout axiom is \cite[lemma 2.13]{Vo}. Let us see the cylinder axiom. If $X\in\simp\mc{C}$, then the natural map $s_0^X:X\boxtimes\Dl[1]\rightarrow X $ is a simplicial homotopy equivalence.
Then, the codiagonal $X\sqcup X\rightarrow X$ factors as the cofibration $X\sqcup X\rightarrow X\boxtimes\Dl[1]$ followed by $s_0^X:X\boxtimes \Dl[1]\rightarrow X\in\mc{W}$.\\
Assume now that $\mc{C}$ is pointed. As seen above, $X\boxtimes\Dl[1]$ is a Brown cylinder for $X$. Therefore, our suspension functor $(\simp\mc{C})[\mc{W}^{-1}]\rightarrow (\simp\mc{C})[\mc{W}^{-1}]$ induced by $\Lambda:\simp\mc{C}\rightarrow \simp\mc{C}$  agrees with the Brown suspension of \cite[theorem 3]{Br}.  Consider a termwise coprojection $i:X\rightarrow Y$ in $\simp\mc{C}$. The Brown cofibration sequence induced by $i$
is
\begin{equation}\label{BrCof}\left( X\stackrel{i}{\rightarrow} Y{\rightarrow} Z=Y\cup_X \ast , \mathtt{b}:Z\rightarrow Z\sqcup \Lambda(X) \right)\end{equation}
where $\mathtt{b}$ is the coaction of \cite[proposition 3]{Br}. By lemma \ref{Pushouts} \textit{i}, we have a natural map $p:\mathtt{Cone}(i)\rightarrow Z$ which is in $\mc{W}$. Moreover, using the description of $\mathtt{b}$ given in \cite[p. 432]{Br}, it is not hard to see that $p$ is compatible with $\mathtt{a}_i$ and $\mathtt{b}$. Consequently, the Brown cofibration sequence (\ref{BrCof}) is isomorphic to the simplicial distinguished triangle induced by
$i$. Conversely, the simplicial distinguished triangle given by $F:X\rightarrow Y$ is isomorphic to the one given by $F':X\rightarrow Y'$, where $F'$ is a cofibration appearing in a factorization $F=t\comp F'$, $t\in\mc{W}$. But, again, this simplicial distinguished triangle is isomorphic to the Brown cofibration sequence of $F'$, so we are done.
\end{proof}

\begin{cor}\label{Triang}\mbox{} Under the hypothesis of previous theorem, the following properties hold.\\
\textbf{\textit{i.}} \emph{(\ref{CogrObj})} makes $\Lambda(X)$ a cogroup object in $(\simp\mc{C})[\mc{W}^{-1}]$, and \emph{(\ref{Coaction})} defines a coaction of
$\Lambda(X)$ on $\mathtt{Cone}(F)$.\\[0.1cm]
\textbf{\textit{ii.}} Simplicial distinguished triangles in $(\simp\mc{C})[\mc{W}^{-1}]$ satisfy the usual `non-stable' axioms for triangulated categories $($see \emph{\cite{Br}}
or \emph{\cite{Q})}.\\[0.1cm]
\textbf{\textit{iii.}}  If $\Lambda:(\simp\mc{C})[\mc{W}^{-1}] \rightarrow (\simp\mc{C})[\mc{W}^{-1}]$ is an equivalence of categories then
$(\simp\mc{C})[\mc{W}^{-1}]$ is a Verdier triangulated category. In particular $(\simp\mc{C})[\mc{W}^{-1}]$ is additive.
\end{cor}

Combining propositions \ref{DlClosed} and \ref{VovWald} we deduce that previous corollary holds for $(\simp\mc{D},\mc{S})$, where $(\mc{D},\mrm{E})$ is a simplicial descent category and $\mc{S}=\mbf{s}^{-1}\mrm{E}$ for some simple functor $\mbf{s}$.

%%%%%%%%%%%%%%%%%%%%%%%%%%%%%%%%%%%%%%%%%%%%%%%%%%%%%%%%%%%%%%%%%%%%%%%%%%%%%%%%%%%%%%%%%%%%%%%%%%%%%%%%%%%%%%%%%%%%%%%%%%%%%%%%%%
%%%%%%%%%%%%%%%%%%%%%%%%%%%%%%%%%%%%%%%%%%%%%%%%%%%%%%%%%%%%%%%%%%%%%%%%%%%%%%%%%%%%%%%%%%%%%%%%%%%%%%%%%%%%%%%%%%%%%%%%%%%%%%%%%%
%%%%%%%%%%%%%%%%%                   CAPITULO
%%%%%%%%%%%%%%%%%%%%%%%%%%%%%%%%%%%%%%%%%%%%%%%%%%%%%%%%%%%%%%%%%%%%%%%%%%%%%%%%%%%%%%%%%%%%%%%%%%%%%%%%%%%%%%%%%%%%%%%%%%%%%%%%%%

\section{The simple functor as homotopical equivalence.}

The results of previous section reveals that the pair $(\simp\mc{D},\mc{S})$ associated with a simplicial descent category $(\mc{D},\mrm{E})$ supports good homotopic properties.
But our aim is to work on $\mc{D}$, not on $\simp\mc{D}$. The key result that makes possible to transfer structure from $\simp\mc{D}$ to $\mc{D}$ is the following.

\begin{thm}\label{equivCat}\mbox{}\\
\textbf{\textit{i.}} The simple functor $\mbf{s}:(\simp\mc{D})[(\simp\mrm{E})^{-1}]\rightarrow \mc{D}[\mrm{E}^{-1}]$ is left adjoint to
$c:\mc{D}[\mrm{E}^{-1}]\rightarrow (\simp\mc{D})[(\simp\mrm{E})^{-1}]$.\\
\textbf{\textit{ii.}} The pair $\mbf{s}:\simp\mc{D}\rightleftarrows \mc{D}:c$ is a homotopical equivalence between
$(\simp\mc{D},\mc{S})$ and $(\mc{D},\mrm{E})$. In particular, $\mbf{s}:\simp\mc{D}[\mc{S}^{-1}]\rightarrow \mc{D}[\mrm{E}^{-1}]$ is an equivalence
of categories.
\end{thm}

\begin{proof}
By (S3) there is a zigzag of natural weak equivalences $\lambda : \mbf{s}\comp c\dashrightarrow 1_{\mc{D}}$ which is then an isomorphism in $Fun(\mc{D},\mc{D})[\mrm{E}^{-1}]$.
To see \textit{i}, it suffices to give $\Phi:1_{\simp\mc{D}}\rightarrow c\comp\mbf{s}$ in $Fun\left(\simp\mc{D},\simp\mc{D}\right)[(\simp\mrm{E})^{-1}]$ such that
\begin{eqnarray}
& \xymatrix@M=4pt@H=4pt@C=25pt{c\ar[r]^-{\Phi_{c}} & c\comp\mbf{s}\comp c\ar[r]^-{c(\lambda)} & c} & \mbox{ is the identity in }Fun(\mc{D},\simp\mc{D})[(\simp\mrm{E})^{-1}]\label{ad1}\\
& \xymatrix@M=4pt@H=4pt@C=25pt{\mbf{s}\ar[r]^-{\mbf{s}(\Phi)} & \mbf{s}\comp c\comp\mbf{s}\ar[r]^-{\lambda_{\mbf{s}}} & \mbf{s}} & \mbox{ is the identity in }Fun(\simp\mc{D},\mc{D})[\mrm{E}^{-1}]\label{ad2}
\end{eqnarray}
Given $X\in\simp\mc{D}$, consider the `total decalage' object associated with $X$ (see \cite[p.7]{I2}). It is a diagram of shape
\begin{equation}\label{DecX}\xymatrix@M=4pt@H=4pt@C=25pt{                                   &                                                                                             &                                                                                                                                              &                                                                                                                                                          &                                                                                                                             & \\
                                      X_{2}\ar@{}[u]|{\vdots}\ar@<0ex>[d]\ar@<1ex>[d] \ar@<-1ex>[d]      & {X_{3}}\ar[l]\ar@<0ex>[d]\ar@<1ex>[d] \ar@<-1ex>[d]\ar@{}[u]|{\vdots}\ar@/_0.75pc/[r]     & {X_{4}} \ar@{}[u]|{\vdots}\ar@<0.5ex>[l] \ar@<-0.5ex>[l]\ar@<0ex>[d]\ar@<1ex>[d] \ar@<-1ex>[d]\ar@/_1pc/[r]\ar@/_0.75pc/[r]                &X_{5} \ar@{}[u]|{\vdots}\ar@<0ex>[l]\ar@<1ex>[l] \ar@<-1ex>[l]\ar@<0ex>[d]\ar@<1ex>[d] \ar@<-1ex>[d]\ar@/_1pc/[r]\ar@/_0.75pc/[r]\ar@{-}@/_1.25pc/[r]   & X_{6}\ar@{}[u]|{\vdots}\ar@<0.33ex>[l]\ar@<-0.33ex>[l]\ar@<1ex>[l]\ar@<-1ex>[l]\ar@<0ex>[d]\ar@<1ex>[d] \ar@<-1ex>[d]     & \cdots \\
                                      X_{1}\ar@<0.5ex>[d] \ar@<-0.5ex>[d] \ar@/^1pc/[u]\ar@/^0.75pc/[u]  & {X_{2}}\ar[l]\ar@<0.5ex>[d] \ar@<-0.5ex>[d] \ar@/^1pc/[u]\ar@/^0.75pc/[u]\ar@/_0.75pc/[r] &  X_{3} \ar@/^1pc/[u]\ar@/^0.75pc/[u]\ar@<0.5ex>[l]\ar@<-0.5ex>[l]\ar@<0.5ex>[d] \ar@<-0.5ex>[d]\ar@/_1pc/[r]\ar@/_0.75pc/[r]               &X_{4} \ar@/^1pc/[u]\ar@/^0.75pc/[u] \ar@<0ex>[l]\ar@<1ex>[l] \ar@<-1ex>[l] \ar@<0.5ex>[d] \ar@<-0.5ex>[d]\ar@/_1pc/[r]\ar@/_0.75pc/[r]\ar@/_1.25pc/[r]  & X_{5}\ar@<0.33ex>[l]\ar@<-0.33ex>[l]\ar@<1ex>[l]\ar@<-1ex>[l] \ar@<0.5ex>[d] \ar@<-0.5ex>[d]\ar@/^1pc/[u]\ar@/^0.75pc/[u] & \cdots \\
                                      X_{0}  \ar@/^0.75pc/[u]                                            &  X_{1}\ar[d] \ar[l]\ar@/_0.75pc/[r] \ar@/^0.75pc/[u]                                      & {X_{2}} \ar[d] \ar@<0.5ex>[l] \ar@<-0.5ex>[l]  \ar@/^0.75pc/[u]  \ar@/_1pc/[r]\ar@/_0.75pc/[r]                                             &X_{3} \ar[d] \ar@<0ex>[l]\ar@<1ex>[l] \ar@<-1ex>[l]\ar@/^0.75pc/[u]\ar@/_1pc/[r]\ar@/_0.75pc/[r]\ar@/_1.25pc/[r]                                        & X_{4}\ar[d] \ar@<0.33ex>[l]\ar@<-0.33ex>[l]\ar@<1ex>[l]\ar@<-1ex>[l] \ar@/^0.75pc/[u]                                     & \cdots \\
                                                                                                         &  X_0 \ar@/_0.75pc/[r]                                                                        & X_1  \ar@<0.5ex>[l]\ar@<-0.5ex>[l]\ar@/_1pc/[r]\ar@/_0.75pc/[r]                                                                         &X_2 \ar@<0ex>[l]\ar@<1ex>[l] \ar@<-1ex>[l]\ar@/_1pc/[r]\ar@/_0.75pc/[r]\ar@/_1.25pc/[r]                                                                   & X_3\ar@<0.33ex>[l]\ar@<-0.33ex>[l]\ar@<1ex>[l]\ar@<-1ex>[l]                                                                 & \cdots }
\end{equation}
\mbox{}\\[0.2cm]
where the morphisms are defined as follows. The $i$-th row is obtained from $X$ by forgetting the last $i+1$ face and degeneracy maps
\begin{equation}\label{decXrow}\xymatrix@1{ X_{i}{\,} & {\;} X_{i+1} \; \ar[l]_{d_0} \ar@/_2pc/[rr]_{s_0} && {\;} X_{i+2} \; \ar@<0.5ex>[ll]^-{d_0} \ar@<-0.5ex>[ll]_-{d_{1}} \ar@/_2.5pc/[rr]_{s_0} \ar@/_1.5pc/[rr]_{s_{1}} && \; X_{i+3} \; \ar@<0ex>[ll]_{d_{1}} \ar@<2ex>[ll]_-{d_{0}} \ar@<-2.25ex>[ll]_-{d_{2}} \ar@/_1.5pc/[rr] \ar@/_2pc/[rr] \ar@/_1pc/[rr] && {\;} X_{i+4}\ar@<0.33ex>[ll]\ar@<-0.33ex>[ll]\ar@<1ex>[ll]\ar@<-1ex>[ll]&\cdots\cdots }\ ,\end{equation}
while the $i$-th column is obtained from $X$ by forgetting the first $i+1$ face and degeneracy maps
\begin{equation}\label{decXcolumn}\xymatrix@1{ X_{i}{\,} & {\;} X_{i+1} \; \ar[l]_{d_{i+1}} \ar@/_2pc/[rr]_{s_{i+1}} && {\;} X_{i+2} \; \ar@<0.5ex>[ll]^-{d_{i+1}} \ar@<-0.5ex>[ll]_-{d_{i+2}} \ar@/_2.5pc/[rr]_{s_{i+1}} \ar@/_1.5pc/[rr]_{s_{i+2}} && \; X_{i+3} \; \ar@<0ex>[ll]_{\small d_{i+2}} \ar@<2ex>[ll]_-{d_{\small i+1}} \ar@<-2.25ex>[ll]_-{\small d_{i+3}} \ar@/_1.5pc/[rr] \ar@/_2pc/[rr] \ar@/_1pc/[rr] && {\;} X_{i+4}\ar@<0.33ex>[ll]\ar@<-0.33ex>[ll]\ar@<1ex>[ll]\ar@<-1ex>[ll]&\cdots\cdots }\end{equation}
Both augmentations have an extra degeneracy: $s_k :X_{i+k}\rightarrow X_{i+k+1}$ for (\ref{decXrow}), and $s_{i}:X_{i+k}\rightarrow X_{i+k+1}$ for (\ref{decXcolumn}).

Let $X\times \Dl$ and $\Dl\times X$ be the bisimplicial objects given by $(X\times \Dl)_{n,m} = X_n$ and $(\Dl\times X)_{n,m} = X_m$. Now, we see diagram (\ref{DecX}) as the bisimplicial object $dec(X)$ given by $dec(X)_{n,m}=X_{n+m+1}$, together with two augmentations $\alpha: dec(X)\rightarrow X\times \Dl$ and $\beta : dec(X)\rightarrow \Dl\times X$.\\
We claim that $(\simp\mbf{s})(\alpha)$ is in $\simp\mrm{E}$. Indeed, given $i\geq 0$, by definition $\alpha_{i,\cdot}$ is the map induced by
the augmentation (\ref{decXcolumn}). Since it has an extra degeneracy, it follows from proposition \ref{ExtraDeg} that $s(m\rightarrow \alpha_{i,m})\in\mrm{E}$ for all $i\geq 0$.
Therefore $(\simp \mbf{s})(\alpha)\in\simp\mrm{E}$.\\
By definition, $(\simp\mbf{s})(X\times\Dl)_n=\mbf{s}(m\rightarrow X_n)=\mbf{s}(c X_n)$.
Axiom (S3) provides a zigzag of natural degreewise weak equivalences $\lambda_{X}:(\simp\mbf{s})(X\times\Dl)\dashrightarrow X$.
On the other hand, $(\simp\mbf{s})(\Dl\times X)$ is just $c\mbf{s}(X)$. We define $\Phi_X :X\dashrightarrow c\mbf{s}(X)$ as
\begin{equation}\label{Phi}
\xymatrix@M=4pt@H=4pt@C=29pt{X   & (\simp\mbf{s})(X\times \Dl) \ar[l]_-{\lambda_X} & (\simp \mbf{s})(dec(X))  \ar[l]_-{(\simp \mbf{s})(\alpha)}
 \ar[r]^-{(\simp \mbf{s})(\beta)} & (\simp\mbf{s})(\Dl\times X)=c\mbf{s}(X)  }
\end{equation}
Given $A$ in $\mc{D}$, $dec(c(A))$ is the constant bisimplicial object equal to $A$, so $\alpha=\beta= 1$. It follows that composition (\ref{ad1}) is equal to the identity.\\
Given $X\in\simp\mc{D}$, (\ref{ad2}) is the top row of the following diagram

$$\xymatrix@M=4pt@H=4pt@C=28pt{
 \mbf{s}(X) & \mbf{s}(\simp\mbf{s})(X\times \Dl) \ar[l]_-{\mbf{s}(\lambda_X)}  & \mbf{s}(\simp \mbf{s})(dec(X)) \ar[l]_-{\mbf{s}(\simp \mbf{s})(\alpha)}
 \ar[r]^-{\mbf{s}(\simp \mbf{s})(\beta)} & \mbf{s}(\simp\mbf{s})(\Dl\times X)=\mbf{s}c\mbf{s}(X) \ar[r]^-{\lambda_{\mbf{s}(X)}} & \mbf{s}(X) \\
 {}  & \mbf{s}(X) \ar[lu]^-{1} \ar[u]_-{\mu_{X\times \Dl}} & \mbf{s}\mrm{D}(dec(X)) \ar[u]_-{\mu_{dec(X)}} \ar[r]^-{\mbf{s}\mrm{D}(\beta)}
\ar[l]_-{\mbf{s}\mrm{D}(\alpha)}  &  \mbf{s}(X) \ar[ru]_-{1} \ar[u]_-{\mu_{\Dl\times X}} & {}  }$$

This is a commutative diagram in $\mc{D}[\mrm{E}^{-1}]$ by (\ref{compatibLambdaMuEquac}).  Therefore, (\ref{ad2}) is the identity if and only if
$\mbf{s}\mrm{D}(\beta)=\mbf{s}\mrm{D}(\alpha)$ in $\mc{D}[\mrm{E}^{-1}]$. But this equality follows from proposition \ref{ExtraDeg}, since by  \cite[proposition 1.6.2]{I2}
the maps $\mrm{D}(\alpha)$ and $\mrm{D}(\beta)$ are simplicially homotopic in a natural way. Therefore, $(\mbf{s},c)$ is an adjoint pair between $(\simp\mc{D})[(\simp\mrm{E})^{-1}]$ and $\mc{D}[\mrm{E}^{-1}]$.\\
Let us see the second statement. By lemma \ref{Wclass} \textit{i},  $\simp\mrm{E}\subset\mc{S}$, so
$c(\mrm{E})\subset \simp\mrm{E}\subset\mc{S}$, while $\mbf{s}(\mc{S})\subset \mrm{E}$ holds by definition.
By (S3), we have a zigzag of natural equivalences $\mbf{s}\comp c \dashrightarrow 1_{\mc{D}}$.
On the other hand, since all arrows in (\ref{Phi}) are in $\mc{S}$ then $\Phi$ provides a zigzag $\Phi^{-1}:c\comp \mbf{s} \dashrightarrow 1_{\simp\mc{D}}$
of natural maps in $\mc{S}$. The fact $(\simp\mbf{s})(\beta)\in\mc{S}$ may be deduced from the above commutative diagram, or as follows. We have that $(\mbf{s}\simp)(\beta)\in\simp\mrm{E} $ because $\beta_{\cdot,i}$ has an
extra degeneracy for each $i\geq 0$. Then $\beta\in (\mbf{s}\simp)^{-1}\mc{S}=(\simp\mbf{s})^{-1}\mc{S}$ by lemma \ref{Wclass} part \textit{ii}.
\end{proof}

\begin{cor}\label{unicSimple}
Let $(\mc{D},\mrm{E})$ be a relative category closed by finite coproducts. Then, all possible simplicial descent structures $(\mbf{s},\mu,\lambda)$ on
$(\mc{D},\mrm{E})$ are unique up to unique isomorphism of $Fun(\simp\mc{D},\mc{D})[\mrm{E}^{-1}]$.
More concretely, given two simplicial descent structures $(\mbf{s},\mu,\lambda)$ and $(\mbf{s}',\mu',\lambda')$ on $(\mc{D},\mrm{E})$, there exists a unique
zigzag of natural weak equivalences $\mbf{s}\dashrightarrow \mbf{s}'$ compatible with $(\mu,\lambda)$ and $(\mu',\lambda')$.
\end{cor}

\begin{proof} First of all, we have that $\mbf{s}\cong\mbf{s}':(\simp\mc{D})[(\simp\mrm{E})^{-1}]\rightarrow\mc{D}[\mrm{E}^{-1}]$ because
they share $c:\mc{D}[\mrm{E}^{-1}]\rightarrow(\simp\mc{D})[(\simp\mrm{E})^{-1}]$ as common right adjoint. Since $\mrm{E}$ is saturated, we deduce
that $\mc{S}=\mbf{s}^{-1}\mrm{E}=\mbf{s}'^{-1}\mrm{E}$. Therefore, as we have seen before, there are zigzags
$\Phi:1_{\simp\mc{D}}\dashrightarrow c\comp \mbf{s}$ and $\Phi':1_{\simp\mc{D}}\dashrightarrow c\comp \mbf{s}'$ of natural maps in $\mc{S}$
such that $(\Phi,\lambda)$ and $(\Phi',\lambda')$ satisfy (\ref{ad1}) and (\ref{ad2}).\\
Then, $\psi=\lambda_{\mbf{s}'}\,\mbf{s}(\Phi'):\mbf{s}\dashrightarrow\mbf{s}'$ is a zigzag of natural weak equivalences compatible with $\lambda$ and $\lambda'$. By (\ref{compatibLambdaMu}), this implies that $\psi$ is also compatible with $\mu_{\Dl\times -}$ and $\mu'_{\Dl\times -}$.
We claim that $\mu$, $\mu'$ are determined by $\mu_{\Dl\times -}$, $\mu'_{\Dl\times -}$. In this case we would deduce that $\psi$ is compatible with $\mu$ and $\mu'$. To see the claim, observe that by proposition \ref{DlClosed} $(\simp\mc{D}, \mc{S})$ is a simplicial descent category with simplicial descent structure $(\mrm{D} , \mu=id,\lambda=id)$, where $\mrm{D}:\simp\simp\mc{D}\rightarrow \simp\mc{D}$ is the diagonal. Hence there is a zigzag of natural maps in $(\mbf{s}\mrm{D})^{-1}\mrm{E}=(\mbf{s}\mbf{s})^{-1}\mrm{E}$ connecting $Z$ to $\Dl\times \mrm{D}(Z)$. Thus $\mu \sim \mu_{\Dl\times\mrm{D}(-)}$ in $Fun(\simp\simp\mc{D},\mc{D})[\mrm{E}^{-1}]$. Analogously $\mu' \sim \mu'_{\Dl\times\mrm{D}(-)}$.\\
To finish, it remains to see that $\psi:\mbf{s}\dashrightarrow\mbf{s}'$ is unique. We must have $\psi_c=\lambda'^{-1}\,\lambda$. Since
$\psi$ is natural, then $\mbf{s}'(\Phi)\,\psi = \psi_{c\, \mbf{s}}\, \mbf{s}(\Phi)$, so
$\psi =\mbf{s}'(\Phi)^{-1} \psi_{c\, \mbf{s}}\, \mbf{s}(\Phi) = \mbf{s}'(\Phi)^{-1} \lambda'^{-1}_{\mbf{s}}\,\lambda_{\mbf{s}}\, \mbf{s}(\Phi)$.
\end{proof}

Combining previous theorem with propositions \ref{DlClosed} \textit{ii} and \ref{VovWald} \textit{i} we deduce the

\begin{cor}
A simplicial descent category is always homotopically equivalent to a Brown category of cofibrant objects.
\end{cor}

\begin{obs}\label{DnoB} Note that not every simplicial descent category $(\mc{D},\mrm{E})$ is itself a Brown category of cofibrant objects.
The reason is that simplicial descent structures are closed by homotopical equivalence (see proposition \ref{EquivHomot}), while Brown
structures of cofibrant objects are not.\\
Consider a simplicial model category $(\mc{M},\mc{W})$ with functorial cofibrant replacements, and such that $\mc{W}$ is closed by finite coproducts.
Its subcategory of cofibrant objects $(\mc{M}_c,\mc{W})$ is then a homotopically equivalent Brown category of cofibrant objects.
By corollary \ref{ModelosCompl} $(\mc{M},\mc{W})$ inherits a simplicial descent structure from $(\mc{M}_c,\mc{W})$.
But $(\mc{M},\mc{W})$ is not necessarily a Brown category of cofibrant objects with the cofibrantions of the model structure, because
it may happen that not all objects in $\mc{M}$ are cofibrant. This is the case, for instance, of  commutative  differential graded algebras over a characteristic 0 field (see \cite{BG}). In the dual setting, one may consider the (Quillen-dual equivalent) model category $(\simp Set,\mc{W})$ where not all objects are fibrant.\\
Note that we obtain in this case a second way to associate with $(\mc{M},\mc{W})$ a homotopically equivalent Brown category of cofibrant objects, namely
$(\simp\mc{M},\mbf{s}^{-1}\mc{W})$.
\end{obs}

\begin{defi} Assume that $(\mc{D},\mrm{E})$ is a pointed simplicial descent category, that is, a simplicial descent category where the initial object $0$ in $\mc{D}$ is also final.
We denote it by $\ast$.\\
The \textit{cone} functor $c:Maps(\mc{D})\rightarrow
\mc{D}$ is defined as
$\mathtt{cone}(f)=\mbf{s}\mathtt{Cone}(c(f))$. Here $\mathtt{Cone}(c(f))$ is the simplicial cone of the simplicial constant map given by $f$ (see definition \ref{DefiCono}).
It is equipped with the natural map $B\rightarrow \mathtt{cone}(f)$ in $\mc{D}[\mrm{E}^{-1}]$ given by $B\stackrel{\lambda_B^{-1}}{\longrightarrow} \mbf{s}c(B)\longrightarrow \mbf{s}\mathtt{Cone}(c(f))$.\\[0.1cm]
The \textit{suspension} functor
${\Sigma}:\mc{D}\rightarrow\mc{D}$ is defined as the simple of the simplicial suspension of $c(A)$, that is,
${\Sigma}(A)=\mbf{s}\Lambda (c(A))=\mathtt{cone}(A\rightarrow\ast)$.\\[0.1cm]
The map $\mathtt{m}_f:\mathtt{cone}(f)\rightarrow \mathtt{cone}(f)\sqcup \Sigma(A)$ in  $\mc{D}[\mrm{E}^{-1}]$ is
$\mathtt{m}_f = \mbf{s} ({\mathtt{a}}_{c(f)})$, where
$$\mathtt{a}_{c(f)}:\mathtt{Cone}(c(f))\rightarrow \mathtt{Cone}(c(f))\sqcup \Lambda(c(A))$$
is given in (\ref{Coaction}).\\[0.1cm]
A distinguished triangle in $\mc{D}[\mrm{E}^{-1}]$ is a pair
$(X\rightarrow Y \rightarrow Z\ ,\ Z\rightarrow Z\sqcup \Sigma X)$
which is isomorphic in $\mc{D}[\mrm{E}^{-1}]$ to a pair of the form
\begin{equation}\label{DTD}\left(A\stackrel{f}{\rightarrow} B \rightarrow \mathtt{cone}(f)\, \ ,  \, \ \mathtt{m}_f:\mathtt{cone}(f)\rightarrow \mathtt{cone}(f)\sqcup \Sigma(A)\right)\end{equation}
\end{defi}

\begin{cor}\label{EstructuraTriangulada} \mbox{}\\[0.1cm]
\textbf{\textit{i.}} Given $A\in\mc{D}$, denote by $\ast_A:A\rightarrow\ast$ the trivial map. Then $\mathtt{m}_{\ast_A}:\Sigma(A)\rightarrow \Sigma(A)\sqcup \Sigma(A)$
makes $\Sigma(A)$ a cogroup object in $\mc{D}[\mrm{E}^{-1}]$. If $f:A\rightarrow B$ then $\mathtt{m}_f$ defines a coaction of
$\Sigma(A)$ on $\mathtt{cone}(f)$.\\[0.1cm]
\textbf{\textit{ii.}} Distinguished triangles in $\mc{D}[\mrm{E}^{-1}]$ satisfy the usual `non-stable' axioms for triangulated categories $($see \emph{\cite{Br}}
or \emph{\cite{Q})}.\\[0.1cm]
\textbf{\textit{iii.}} Distinguished triangles \emph{(\ref{DTD})} are natural with respect to diagram categories. In other words, given a small category $I$, there is a natural notion of distinguished triangles in $(I\mc{D})[\mrm{E}^{-1}]$. In addition, if $\phi:I\rightarrow J$ is a functor of small categories, then
$\phi^\ast :(J\mc{D})[\mrm{E}^{-1}]\rightarrow (I\mc{D})[\mrm{E}^{-1}]$ preserves distinguished triangles.\\[0.1cm]
\textbf{\textit{iv.}}  Assume moreover that $\Sigma:\mc{D}[\mrm{E}^{-1}] \rightarrow \mc{D}[\mrm{E}^{-1}]$ is an equivalence of categories. Then
$\mc{D}[\mrm{E}^{-1}]$ is a Verdier triangulated category. In particular $\mc{D}[\mrm{E}^{-1}]$ is additive.\\[0.1cm]
\textbf{\textit{v.}} If $\Sigma$ is an isomorphism in $Fun(\mc{D},\mc{D})[\mrm{E}^{-1}]$, then $(I\mc{D})[\mrm{E}^{-1}]$ is a Verdier triangulated category for each small
category $I$, and each $\phi:I\rightarrow J$ induces a triangulated functor $\phi^\ast :(J\mc{D})[\mrm{E}^{-1}]\rightarrow (I\mc{D})[\mrm{E}^{-1}]$.
\end{cor}

\begin{proof}
Parts \emph{i}, \emph{ii} and \emph{iv} are obtained combining proposition \ref{DlClosed}, theorem \ref{equivCat} and corollary \ref{Triang}. Parts \emph{iii} and \emph{v} are consequences of the previous ones and the fact that simplicial descent structures are inherited by diagram categories by proposition \ref{DescensoFuntores}.
\end{proof}

%%%%%%%%%%%%%%%%%%%%%%%%%%%%%%%%%%%%%%%%%%%%%%%%%%%%%%%%%%%%%%%%%%%%%%%%%%%%%%%%%%%%%%%%%%%%%%%%%%%%%%%%%%%%%%%%%%%%%%%%%%%%%%%%%%
%%%%%%%%%%%%%%%%%%%%%%%%%%%%%%%%%%%%%%%%%%%%%%%%%%%%%%%%%%%%%%%%%%%%%%%%%%%%%%%%%%%%%%%%%%%%%%%%%%%%%%%%%%%%%%%%%%%%%%%%%%%%%%%%%%
%%%%%%%%%%%%%%%%%                   CAPITULO
%%%%%%%%%%%%%%%%%%%%%%%%%%%%%%%%%%%%%%%%%%%%%%%%%%%%%%%%%%%%%%%%%%%%%%%%%%%%%%%%%%%%%%%%%%%%%%%%%%%%%%%%%%%%%%%%%%%%%%%%%%%%%%%%%%


\begin{thebibliography}{99}
\begin{small}
%
\bibitem[B]{B} A. A. Beilinson, {\it Notes on absolute Hodge cohomology.  Applications of algebraic $K$-theory to algebraic geometry and number theory}, Contemp. Math., \textbf{55},
Amer. Math. Soc., (1986) p. 35-68.
%
\vspace{-0.1cm}
%
\bibitem[BG]{BG} A. K. Bousfield and A. M. Gugenheim, \textit{On PL de Rham theory and rational homotopy type}, Mem. Amer. Math. Soc. \textbf{179} (1976).
%
\vspace{-0.1cm}
%
\bibitem[Br]{Br} K. S. Brown, {\it Abstract homotopy theory and generalized sheaf cohomology}, Trans. Amer. Math. Soc. \textbf{186} (1974), p. 419-458.
%
\vspace{-0.1cm}
%
\bibitem[CN]{CN} C. Casacubertas and A. Neeman, {\it Brown representability does not come for free}, Math. Res. Lett. \textbf{16}  (2009),  no. 1, p. 1-5.
%
\vspace{-0.1cm}
%
\bibitem[DeII]{DeII} P. Deligne, \textit{Th\'{e}orie de Hodge II}, Publ. Math. I.H.E.S.,
\textbf{40} (1971), p. 5-57.
%
\vspace{-0.1cm}
%
\bibitem[DeIII]{DeIII} P. Deligne, \textit{Th\'{e}orie de Hodge III}, Publ. Math. I.H.E.S.,
\textbf{44} (1975), p. 2-77.
%
\vspace{-0.1cm}
%
\bibitem[DP]{DP} A. Dold and D. Puppe, {\it Homologie nicht-additiver funktoren. Anwendungen}, Ann. Inst. Fourier, \textbf{11} (1961), p. 201-312.
%
\vspace{-0.1cm}
%
\bibitem[DHKS]{DHKS} W. G. Dwyer, P. S. Hirschhorn, D. M. Kan, and J. H. Smith, \textit{Homotopy limit functors
on model categories and homotopical categories}, Math. Surveys and Monographs, \textbf{113}, AMS, Providence, (2004).

\vspace{-0.1cm}
%
\bibitem[EM]{EM} S. Eilenberg and J.C. Moore, \textit{Homology and fibrations. I. Coalgebras, cotensor product and its derived functors},
Comment. Math. Helv. \textbf{40} (1966), p. 199-236.
%
\vspace{-0.1cm}
%
\bibitem[F]{F} P. Freyd, \textit{Abelian categories}, Harper and Row, New York (1964).
%
\vspace{-0.1cm}
%
\bibitem[GJ]{GJ} P.G. Goers y J.F. Jardine, \textit{Simplicial homotopy theory}, Birkh{\"a}user Verlag, Basel, 1999.
%
\vspace{-0.1cm}
%
\bibitem[GN]{GN} F. Guill\'{e}n and V. Navarro Aznar, \textit{Un crit\`{e}re d'extension des foncteurs d\'{e}finis sur les sch\'{e}mas lisses}, Publ. Math. I.H.E.S., \textbf{95} (2002), p. 1-91.
%
\vspace{-0.1cm}
%
\bibitem[H]{H} P. S. Hirschhorn, \textit{Model Categories and their localizations}, Math. Surveys and Monographs, \textbf{99}, AMS, Providence (2002).
%
\vspace{-0.1cm}
%
\bibitem[Hb]{Hb} A. Huber, \textit{Mixed motives and their realization in derived categories}, Lect. Notes in Math., \textbf{1604}, Springer, Berlin (1995).
%
\vspace{-0.1cm}
%
\bibitem[I1]{I1} L. Illusie, \textit{Complexe cotangent et d\'{e}formations I}, Lect. Notes in Math., \textbf{239}, Springer, Berlin (1971).
%
\vspace{-0.1cm}
%
\bibitem[I2]{I2} L. Illusie, \textit{Complexe cotangent et d{\'e}formations II}, Lect. Notes in Math., \textbf{283}, Springer, Berlin (1972).

\vspace{-0.1cm}
%
\bibitem[ML]{ML} S. MacLane, {\it Categories for the working mathematician}, GTM, \textbf{5}, Springer, Berlin, 1971.
%
\vspace{-0.1cm}
%
\bibitem[May]{May} J. P. May, {\it Simplicial objects in Algebraic Topology}, Van Nostrand, Princeton, 1967.
%
\vspace{-0.1cm}
%
\bibitem[PS]{PS} C. Peters and J. Steenbrink, {\it Mixed Hodge structures}, Ergebnisse der Math.; Ser. of Modern Surveys in Mathematics, \textbf{52}, Springer-Verlag, Berlin  (2008).
%
\vspace{-0.1cm}
%
\bibitem[Q]{Q} D. Quillen, \textit{Homotopical Algebra}, Lect. Notes in Math., \textbf{43}, Springer, Berlin (1967).
%
\vspace{-0.1cm}
%
\bibitem[R1]{R1} B. Rodr\'{\i}guez-Gonz\'{a}lez, \textit{Categor\'ias de descenso simplicial}, PhD thesis, University of Seville, 2007. Translated and revised version at \href{http://arxiv.org/abs/0804.2154}{arXiv:0804.2154v1}.
%
\vspace{-0.1cm}
%
\bibitem[R2]{R2} B. Rodr\'{\i}guez-Gonz\'{a}lez, \textit{Realizable homotopy colimits}, preprint 2011, \href{http://arxiv.org/abs/1104.0646}{arXiv:1104.0646}.
%
\vspace{-0.1cm}
%
\bibitem[Ru]{Ru} L. Rubi{\'o} i Pons, \textit{Model categories and cubical descent}, Bol. Soc. Mat. Mexicana, \textbf{13} (2007), p. 293-305.

\vspace{-0.1cm}
%
\bibitem[T]{T} R. W. Thomason, {\it Algebraic K-theory and etale cohomology},  Ann. Sci. de L'E.N.S \textbf{18}, (1985), p. 437-552.
%
\bibitem[Vo]{Vo} V. Voevodsky, {\it Simplicial radditive functors}, Journal of K-theory, \textbf{5} (2010), p. 201-244.
%
\vspace{-0.1cm}
%
\bibitem[W]{W} C. Weibel, {\it An introduction to homological algebra},
Cambridge studies in advanced math. \textbf{38} (2003), Cambridge
Univ. press.
%
\end{small}
\end{thebibliography}
\end{document}